\documentclass{article}
\usepackage[utf8]{inputenc}
\usepackage{lmodern}
\usepackage[T1]{fontenc}
\usepackage{comment}
\usepackage{lscape}
\usepackage{amsmath}	
\usepackage{amssymb}
\usepackage{amsthm}
\usepackage{cases}
\usepackage{graphicx}
\usepackage{caption}

\usepackage{makecell}

\usepackage{multirow}
\usepackage{amsfonts}
\usepackage{float} 
\usepackage{diagbox}

\usepackage{fullpage}
\usepackage{url}
\usepackage{rotating}
\usepackage{eurosym}
\usepackage{wrapfig}
\usepackage[final]{pdfpages} 
\usepackage{epstopdf} 
\usepackage{subfig} 
\usepackage[a4paper]{geometry}
\usepackage[nottoc]{tocbibind} 
\usepackage{placeins}
\usepackage{float}

\usepackage{listings}
\usepackage{color, colortbl}

\usepackage{hyperref}
\usepackage{xcolor}
\usepackage{titlesec}
\usepackage[ruled]{algorithm2e}

\titleformat{\paragraph}
{\normalfont\normalsize\bfseries}{\theparagraph}{1em}{}
\titlespacing*{\paragraph}
{0pt}{3.25ex plus 1ex minus .2ex}{1.5ex plus .2ex}

\hypersetup{
colorlinks,%
citecolor=black,%
filecolor=black,%
linkcolor=black,%
urlcolor=black
} 

\newcommand{\di}{\textbf{DI01 }}
\newcommand{\COV}{\textbf{\text{cov}}}

\begin{document}
\title{A graph clustering approach to localization for adaptive covariance tuning in data assimilation based on state-observation mapping}

\author{Sibo Cheng$^{1,2}$, Jean-Philippe Argaud$^{1}$, Bertrand Iooss$^{1,3}$,  Angélique Pon{ç}ot$^{1}$, Didier Lucor$^{2}$ \\
\\
        \small $^{1}$ EDF R\&D \\
        \small $^{2}$ LIMSI, CNRS, Univ.Paris-Sud, Université Paris-Saclay\\
        \small $^{3}$ Institut de Mathématiques de Toulouse, Université Paul Sabatier\\
}
\maketitle
\date{}

\begin{abstract}
An original graph clustering approach to efficient localization of error covariances is proposed  within an ensemble-variational data assimilation framework. Here the localization term is very generic and refers to the idea of breaking up a global assimilation into subproblems. This unsupervised localization  technique based on a linearized
state-observation measure is general and does not rely on any prior information such as relevant spatial scales, empirical cut-off radius or homogeneity assumptions. It automatically segregates the state and observation variables in an optimal number of clusters (otherwise named as subspaces or communities), more amenable to scalable data assimilation.
The application of this method does not require underlying block-diagonal structures of prior covariance matrices. In order to deal with inter-cluster connectivity, two alternative data adaptations are proposed. 
Once the localization is completed, an adaptive covariance diagnosis and tuning is performed within each cluster. 
Numerical tests show that this approach is less costly and more flexible than a global covariance tuning, and most often results in more accurate background and observations error covariances. 
\end{abstract}


\section{Introduction}

Much attention has been devoted in data assimilation to the modeling of background or/and observation error covariance representation, $\textbf{B}$ or/and $\textbf{R}$ respectively. These prior covariance matrices can be estimated from the help of a correlation kernel (e.g. \cite{dance2013}, \cite{Poncot2013}) or a diffusion operator (e.g. \cite{Mirouze2010}), or they may be improved by the NMC method (\cite{F.Parrish1992}), ensemble methods (\cite{Clayton2012}), time lagged modelling (\cite{Thomas2017}) or some iterative methods for which a part of $\textbf{B}$ and/or $\textbf{R}$ is supposed to be perfectly known  e.g. \cite{Desroziers01}, \cite{Desroziers2005}, \cite{cheng2019}. These approaches quite often rely on converged state ensemble statistics, noiseless dynamical system or assumption of error amplitude (\cite{cheng2019}). These conditions are sometimes difficult to be satisfied in an industrial context which holds limited accessibility to historical data. 

In fact, in the case of high-dimensional systems, if the ensemble is too small, the estimated full error covariance matrix will most likely be rank deficient and polluted by spurious correlations at long distances. The idea of {\em localization} relies on the intuitive idea that distant states of the system are more likely to be independent, at least for sufficiently short time scales. For many application fields, such as geophysical sciences, the system state depends on spatial coordinates, which are easily monitored and comprehended, and it is then possible to {\em spatially} localize the analysis. For other systems (e.g. the interchannel radiance observation, \cite{Garand2007}), the correlation between different ranges/scales of the state variables may not be directly interpreted in terms of spatial distances and the assumption of weak long-distance correlations might be less relevant (\cite{Buehner2015}). In this paper, we will refer to the more generic ``long-range correlation'' expression instead.  Also, there might be situations for which a prior covariance structure has limited spatial extent, that is smaller than the support of the observation operator that maps state space to observations space.
In this case, non-local observations, i.e. observations that cannot be really allocated to one specific spatial location, because they may result from spatial averages of linear or non-linear functions of the system variables can have a large influence on the assimilation, cf. work of  \cite{Leeuwen2019}.

Existing localization methods are mainly two kinds: -- covariance localization and -- domain localization.
The first family of localization methods is implicit in the sense that it directly works on a regularization of the covariance matrix (e.g.\cite{Yoshida2018}) that is operated using a Schur matrix product with certain short-range predefined correlation matrices, such as the Gaspari-Cohn correlation function (\cite{Gaspari1998}). The Schur product theorem ensures the (semi)definitiveness of the new matrix and therefore avoids the introduction of spurious long-range correlation. These methods have been widely improved, in particular for ensemble-based Kalman filters (EnKF) (e.g. \cite{Farchi2019}) where the covariance localization is crucial to produce more accurate analyses. 

The second class of families (i.e. domain localization) is more explicit and performs data assimilation for each state variable by using only a local subset of available observations, typically within a fixed range of this point, 
as for instance the LETKF method of \cite{Hunt2007} in a EnKF framework.
In this case, a relevant localization length must be carefully chosen. This is the main disadvantage of the approach. For instance, if this length is chosen too small, some important short- to medium-range correlation will be falsely neglected.

Recent works have shown that a local diagnosis/correction of background-observation covariance computation could be helpful for improving the forecast quality of the global system, e.g. \cite{Waller2017}, as well as reducing the computational cost. From the point of view of an observation, it introduces the concepts of –  {\em domain of dependence}, i.e. the set of elements of the model state that are used to predict the model equivalent of this observation; and of – the {\em region of influence}, i.e. the set of analysis states that are updated in the assimilation using this observation. According to \cite{Waller2017}, difficulties appear with the domain localization when the region of influence is far offset from the domain of dependence. In fact, the former may be imposed based on prior assumptions while the later is obtained from the linearized transformation operator, which depicts how the state variables are ``connected'' via the observations. 

To summarize, the main drawback of localization performed that ways is that the cut-off or distance thresholds are somewhat arbitrarily chosen  and may result in removal of true physical long-range correlations, thus inducing imbalance in the analysis as mentioned by \cite{Greybush2011}. In fact, data assimilation often deals with non-uniform inhomogeneous error fields, containing underlying structure due to the heterogeneity of the data, which calls for powerful unsupervised localization schemes. 
This points to the relevance of more efficient and less arbitrary segregation operators. One of the main objectives of  unsupervised learning is to split a data set into two or more classes based on a similarity measure over the data, without resorting to any {\em a priori} information on how the community clustering should be done (see \cite{hastie2001}, section 14).
In this study, we choose to segregate the state variables directly based on the information provided by the state-observation mapping.
This unifying approach avoids potential conflicts between the region of influence and the domain of dependence of the localized assimilation.

A first original idea  of our work, is to turn to efficient localization strategies based on {\em graph clustering} theory, which are able to automatically detect several clusters or ``communities'' (we will also refer to them as ``subspaces'' in the state and observation space) of state variables and corresponding observations. This clustering of variables will allow more local assimilation, likely to be more flexible and accurate  than  a standard global assimilation technique.  
Some elements of graph theory have already been introduced in data assimilation algorithms, for example in \cite{Ihler2005} where each edge represents a conditional dependence among state variables and observations. In this paper, we adopt a different approach where community detection algorithms are applied  in a network of state variables in order to improve the  subspaces segmentation. This network, called an observation-based state network, will only depend on the linearized transformation operator $\textbf{H}$ between state variables and observations. More precisely, our objective is to classify the state variables represented by the same observation to the same subspace, regardless of their distance.

Assuming the graph clustering approach has been efficient for localizing several state communities, data assimilation/covariances diagnosis still has to be performed within each of them. To this end, it is crucial to assign appropriate subset of observations to each community of state variables, and different modeling and computational approaches are possible along those lines.

Similar to \cite{Waller2017}, our strategy is  then to conduct standard error covariance estimation in the detected subspaces. As mentioned previously, remarkable effort has been made on posterior diagnosing and iterative adjustment of error covariance quantification, especially by the meteorology community. 
(e.g. \cite{Desroziers01}, \cite{Desroziers2005}). 
Among these works, the tuning method of Desroziers and Ivanov  (also known as \di), which consists of finding a fixed point for the assimilated state by adjusting the ratio between background and observation covariance matrices amplitude without modifying their structures, is widely adopted in numerical weather prediction. This approach presents the  advantage that it can be implemented in a static data assimilation problem or at any step of a dynamical data assimilation process for both variational methods and  Kalman-type filtering, even with limited background/observation data, which is not the case for a full covariance estimation/diagnosis (e.g. \cite{Desroziers2005}). The later is based on statistics of prior and posterior innovation quantities.
In fact, the deployment  of \di~in subspaces has already been introduced in \cite{Chapnik2006} for block diagonal structures of $\textbf{B}$ and $\textbf{R}$. In this paper, we adopt a \di~approach that we extend to a more general approach, where the block diagonal structure of the covariances matrix  is no longer required, but covariance between extra-diagonal blocks remains accounted for. 

The paper is organized as follows: the standard formulation of data assimilation is introduced, as well as its resolution in the case of a linearized Jacobian matrix, in section \ref{sec:da_principe}
We then explain how this Jacobian matrix, considered as a state-observation mapping, can be used to build an observation-based state network where the subspaces decomposition is carried out by applying graph-based community detection algorithms. The localized version of \di~is then introduced (section \ref{sec:Desroziers Ivanov tuning algorithm}) and investigated in a twin experiments framework (section \ref{sec:Illustration with numerical experiments}). We close the paper with a discussion (section \ref{sec:discussion}).

\section{Data assimilation framework}\label{sec:da_principe}

The goal of data assimilation algorithms is to correct the state $\textbf{x}$ of a a dynamical system with the help of a prior estimation $\textbf{x}_b$ and an observation vector $\textbf{y}$, the former being often provided by expertise or a numerical simulation code. This correction brings the state vector closer to its true value denoted by $\textbf{x}_t$, also known as the true state. In this paper, each state component $\textbf{x}_i$ is called a state variable while $\textbf{y}_j$ is called an observation.  
The principle of data assimilation algorithms is to find an optimally weighted combination of $\textbf{x}_b$ and $\textbf{y}$ by optimizing the minimum cost of a cost function $J$ defined as:
\begin{align}
    J(\textbf{x})&=\frac{1}{2}(\textbf{x}-\textbf{x}_b)^T\textbf{B}^{-1}(\textbf{x}-\textbf{x}_b) \notag \\
    & + \frac{1}{2}(\textbf{y}-\mathcal{H}(\textbf{x}))^T \textbf{R}^{-1} (\textbf{y}-\mathcal{H}(\textbf{x})) \label{eq_3dvar}\\
    & = J_b(\textbf{x}) + J_o(\textbf{x})
\end{align}
 where the observation operator $\mathcal{H}$ denotes the mapping from the state space to the one of observations. $\textbf{B}$ and $\textbf{R}$ are the associated error covariance matrices, i.e.
 \begin{align}
     \textbf{B} & = \COV(\epsilon_b, \epsilon_b),\\
     \textbf{R} & = \COV(\epsilon_y, \epsilon_y),
 \end{align}
 where
  \begin{align}
     \epsilon_b & = \textbf{x}_b - \textbf{x}_t, \\
     \epsilon_y & = \mathcal{H}(\textbf{x}_t)-\textbf{y}.
 \end{align}
 Their inverse matrices, $\textbf{B}^{-1}$ and $\textbf{R}^{-1}$, represent the weights of these two information sources in the objective function. Prior errors, $\epsilon_b$ and $\epsilon_y$, are supposed to be centered Gaussian variables in data assimilation, thus they can be perfectly characterized by the covariance matrices, i.e.
 
 \begin{align}
     \epsilon_b & \sim \mathcal{N} (0, \textbf{B}), \\
     \epsilon_y & \sim \mathcal{N} (0, \textbf{R}).
 \end{align}
 
The two covariance matrices $\textbf{B}$ and $\textbf{R}$, which are difficult to know perfectly \textit{a priori}, play essential roles in data assimilation. The state-observation mapping $\mathcal{H}$  is possibly nonlinear in real applications. However, for the sake of simplicity, a linearization of $\mathcal{H}$ is often required to evaluate the posterior state and its covariance. The linearized operator $\textbf{H}$, often known as the Jacobian matrix in data assimilation, can be seen as a mapping from the state space to the one of observation. 
 
 In the case where $\mathcal{H} = \textbf{H}$ is linear and the covariances matrices $\textbf{B}$ and $\textbf{R}$ are well known, the optimization problem (\ref{eq_3dvar})  can be perfectly solved  by linear formulation of BLUE:
 \begin{equation}
	\textbf{x}_a=\textbf{x}_b+\textbf{K}(\textbf{y}-\textbf{H} \textbf{x}_b) \label{eq:BLUE_1}
\end{equation}
which is also equivalent to a maximum likelihood estimator. 
The Kalman gain matrix $\textbf{K}$ is defined as:
\begin{equation}
	\textbf{K}=\textbf{B} \textbf{H}^T (\textbf{H} \textbf{B} \textbf{H}^T+\textbf{R})^{-1}. 
	\label{eq:Kgain_BLUE}
\end{equation}
 
 Several diagnosis or tuning methods, such as the ones of \cite{Desroziers2005}, \cite{Desroziers01}, \cite{Dreano2017} have been developed to improve the quality of covariance estimation/construction. Much effort has also been devoted to apply these methods in subspaces (e.g. \cite{Waller2017}, \cite{Sandu2015}). The subspaces are often divided by the physical nature of state variables or their spatial distance. The prior estimation errors are often considered as uncorrelated among different subspaces. A significant disadvantage of this approach is that the cut-off correlation radius remains difficult to determine and the hypothesis of no error correlation among distant state variables is not always relevant depending on the application.

\section{State-observation localization based on graph clustering methods}

For the purpose of the simplicity of implementation, representing state variables/observations by block diagonal matrices is sometimes used in data assimilation (for example, see \cite{Chabot2015}). In this case, only uncorrelated state variables can be separated. In this work,
we are interested in applying  covariance diagnosis methods in subspaces identified from the state-observation mapping, 
and we wish to make no assumption of block diagonal structures of the covariance.

The state subspaces will be detected thanks to an unsupervised graph clustering learning technique.
Here, the graph will be formed by a set of vertices (i.e. the state discrete nodes) and a set of edges (based on a similarity measure over the state variables-observations mapping) connecting pairs of vertices. The graph clustering will automatically group the vertices of the graph into clusters taking into consideration the edge structure of the graph in such a way that there should be many edges within each cluster and relatively few between the clusters.

\subsection{State space decomposition via graph clustering algorithms}


\subsubsection{Principles}

Here, the idea is to perform a localization by segregating the state vector $\textbf{x}\in \mathbb{R}^{n_{\textbf{x}}}$ (we drop the background or analysed subscript for the ease of notation) into a partition $\mathcal{C}$ of subvectors: $\mathcal{C}=\{\textbf{x}^1,\textbf{x}^2,\ldots, \textbf{x}^p \}$, each $\textbf{x}^i$ being non-empty. We will call later $\mathcal{C}$ a {\em clustering} and the elements $\textbf{x}^i$ {\em clusters}.
Similarly to the standard localization approach, for each identified subset of state variables, it will then be necessary to identify an associated subset of observations: $\{\textbf{y}^1,\textbf{y}^2,\ldots, \textbf{y}^p \}$.  

In the work of \cite{Waller2016}, a threshold of spatial distance $\tilde{r}$ is {\em arbitrarily} imposed \textit{a priori} to define local subsets of state variables influenced by each observation during the data assimilation updating. In other words, each observation component $\textbf{y}_i$, of the complete vector $\textbf{y}$,
is only supposed to influence the updating of a subset of state variables  within the spatial range of  $\tilde{r}$. This subset of state variables $\mathcal{R}_{\textrm{influence}}(\textbf{y}_i) = \{ \textbf{x}_k: \phi(\textbf{y}_i, \textbf{x}_k) \leq \tilde{r} \}$, where $\phi$ measures some spatial distance,
 is called the {\em region of influence} of $\textbf{y}_i$. 

However, that method faces a significant difficulty when the Jacobian matrix $\textbf{H}$ is dense or non local, i.e. the updating of state variables depends on observations out of the region of influence. 
In fact, the non-locality of matrix $\textbf{H}$ may contain terms that will induce a ``connection'' between state variables and observations beyond the critical spatial range $\tilde{r}$. The {\em domain of dependence} defined as:
\begin{align}
    \mathcal{D}_{\textrm{dependence}}(\textbf{y}_i)= \{ \textbf{x}_k: \textbf{H}_{i,k} \neq 0 \},
\end{align}
is introduced to quantify the range of this state-observation connection which is purely decided by $\textbf{H}$ instead of the spatial distance. \cite{Waller2016} have shown that problems may occur in the covariance diagnosis when $\mathcal{R}_{\textrm{influence}}(\textbf{y}_i) $ and $\mathcal{D}_{\textrm{dependence}}(\textbf{y}_i) $ do not overlap. {This incoherence  not only impacts the assimilation accuracy but also the posterior covariance estimation. This phenomenon is also highlighted and studied in the work of \cite{Leeuwen2019} where the author proposes an extra step to assimilate observations outside the region of influence.}

\subsubsection{Observation-based state connections}
 Rather than considering the region of influence, our proposed approach uses a clustering strategy directly based  on the domain of dependence, i.e. taking advantage of the particular structure of the Jacobian $\textbf{H}$. The main idea is to separate the ensemble of state variables into several subsets regarding their occurrence in the domains of dependence of different observations. 
 We introduce the notion of observation-based connection between two state variables $\textbf{x}_i$ and $\textbf{x}_j$ when they appear in the domain of dependence of the same observation $\textbf{y}_k$, i.e.
\begin{align}
    \exists k,\quad \textrm{such that} \quad \frac{\partial \textbf{y}_k}{\partial \textbf{x}_i} \neq 0, \quad  \quad \frac{\partial \textbf{y}_k}{\partial \textbf{x}_j} \neq 0.
\label{eq:link}
\end{align}
For a linearized state-observation operator $\textbf{H}$, it is simply equivalent to 
\begin{align}
\exists k, \quad \textrm{such that} \quad \textbf{H}_{k,i} \neq 0, \quad \quad \textbf{H}_{k,j} \neq 0.
\label{eq:link-H}
\end{align}

Our goal is to determine if we can group the state variables which are strongly connected based on the observations, regardless of their spatial distance. In order to do so, we define the strength of this connection for each pair of state variables. 
In fact, $\textbf{H}$ can be seen as some form of weighted mapping between the space of state variables and the one of observations, which allows us to define the connection strength, quantified by a function $\mathcal{S}$ as a sum of conjugated coefficients multiplication in $\textbf{H}$, i.e,
\begin{eqnarray}
    \mathcal{S} & : & \mathbb{R}^{n_{\textbf{x}}} \times \mathbb{R}^{n_{\textbf{x}}} \mapsto \mathbb{R}^{+} \\
    & & \mathcal{S}(\textbf{x}_i,\textbf{x}_j) \equiv \mathcal{S}_{i,j} =\sum\limits_{k, i\neq j} |\textbf{H}|_{k,i} |\textbf{H}|_{k,j}, \nonumber
\end{eqnarray}
where $|\cdot|$ represents the absolute value function on the whole matrix, so the function is symmetric. Moreover, we will assume that the function is null when measuring the connection strength of one state variable with itself.

We now consider an undirected  graph $\mathcal{G}$ that is a pair of sets $\mathcal{G}=(\textbf{x},\textbf{E})$, where $\textbf{x}$ plays the role of the set of vertices (the number of vertices $n_{\textbf{x}}$ is the order of the graph) and the set $\textbf{E}$ contains the edges of the graph (the edge cardinality, i.e. $|\textbf{E}|=m$ represents the size of the graph). Each edge is an unordered pair of endpoints $\{\textbf{x}_k,\textbf{x}_l \}$.
We are going to use our measure $\mathcal{S}$ as a weight function to define the weighted version of the graph $\mathcal{G}_{\mathcal{S}}=(\textbf{x},\textbf{E},\mathcal{S})$. This translates into the weighted adjacency matrix $\textbf{A}_{\mathcal{G}_{\mathcal{S}}}$ of the graph, that is a $n_{\textbf{x}}\times n_{\textbf{x}}$ matrix $\displaystyle \textbf{A}_{\mathcal{G}_{\mathcal{S}}}=(a^{\mathcal{G}_{\mathcal{S}}}_{\textbf{x}_i,\textbf{x}_j})$:
 \begin{equation}
a^{\mathcal{G}_{\mathcal{S}}}_{\textbf{x}_i,\textbf{x}_j} =\left\{
\begin{array}{c l}	
     & \mathcal{S}_{i,j} \; \text{if} \; \{\textbf{x}_i,\textbf{x}_j\}\in \textbf{E},\\
     & 0 \quad \textrm{otherwise.}
\end{array}\right.
\end{equation}
This matrix will be useful to perform the graph clustering.

 Each edge of the graph thus represents the connection strength between two state variables. For some problems, it is possible to organize the graph into clusters, with many edges joining vertices of the same cluster and comparatively few edges joining vertices of different clusters.
 We have the partition $\mathcal{C}=\{\textbf{x}^1,\textbf{x}^2,\ldots, \textbf{x}^p \}$ of $\textbf{x}$, and we identify a cluster $\textbf{x}^{i}$ with a node-induced subgraph of $\mathcal{G}_{\mathcal{S}}$, i.e. the subgraph $\mathcal{G}_{\mathcal{S}}\left [\textbf{x}^i \right ] \mathrel{\mathop:}=\left ( \textbf{x}^i,\textbf{E}(\textbf{x}^i),\mathcal{S}_{|\textbf{E}(\textbf{x}^i)} \right )$, where $\textbf{E}(\textbf{x}^i)\mathrel{\mathop:}=\left \{ \left \{ \textbf{x}_k,\textbf{x}_l\right \}\in \textbf{E}: \textbf{x}_k,\textbf{x}_l \in \textbf{x}^i \right \}$. So $\textbf{E}(\mathcal{C}) \mathrel{\mathop:}= \bigcup_{i=1}^p \textbf{E}(\textbf{x}^i)$ is the set of intra-cluster edges and $\textbf{E}\textbackslash \textbf{E}(\mathcal{C})$ is the set of inter-cluster edges of cluster $\textbf{x}^i$ respectively, with $\left | \textbf{E}(\mathcal{C})  \right | = m(\mathcal{C})$ and $\left | \textbf{E}\textbackslash \textbf{E}(\mathcal{C}) \right |=\bar{m}(\mathcal{C})$, while 
 $\textbf{E}(\textbf{x}^i,\textbf{x}^j)$ denotes the set of edges connecting nodes in $\textbf{x}^i$ to nodes in $\textbf{x}^j$. 
It is important to stress that the identification of structural clusters is made easier if graphs are {\em sparse}, i. e. if the number of edges $m$ is of the order of the number of nodes $n_{\textbf{x}}$ of the graph (\cite{FORTUNATO2010}).

\subsubsection{Clustering algorithms}

One of the main paradigms of clustering is to find groups/clusters intra-cluster density vs. inter-cluster sparsity. Despite the fact that many problems related to clustering are $\textbf{NP}$-hard problems, there exist many approximation methods for graph-based community detection, such as the Louvain algorithm (\cite{Blondel2008}) and the Fluid community algorithm (\cite{pares2018}). These methods are mostly based on random walks or centrality measures in a network with the advantage of low computational cost. The use of graph theory in numerical simulation problems such as the Cuthill–McKee algorithm (\cite{Cuthill1969}) already exists, for instance for sorting multidimensional grid points in a more efficient way (in terms of reducing the matrix band). In this paper, we introduce a different approach with the objective of identifying observation-based state variable communities which will be later considered as state subsets in covariance tuning. The community detection is performed on the observation-based state network, regardless of the algorithms chosen. Considering the computational cost, the Fluid community detection algorithm proposed by \cite{pares2018} could be an appropriate choice for sparse transformation matrix because its complexity is linear to the number of edges in the network, i.e. $\mathcal{O}(|\textbf{E}|)$.

In real applications of graph theory, the number of optimal cluster $p$ is often not known in advance. Finding appropriate cluster number remains a popular research topic.
Several methods have been developed in order to propose some objective functions with notion of optimal coverage, performance or inter-cluster conductance, e.g. the Elbow method (\cite{KETCHEN1996}) or the Gap statistic method (\cite{Tibshirani2001}). For instance the following performance metric will be used later for the experiments in section \ref{sec:Experimental results}:
\begin{equation}
\text{performance}\mathrel{\mathop:}= \frac{m(\mathcal{C})+\bar{m}^c(\mathcal{C})}{\frac{1}{2}n_{\textbf{x}}(n_{\textbf{x}}-1)}.
\end{equation}
It represents the fraction of node pairs that are clustered correctly, i.e. those connected node pairs that are in the same cluster and those non-connected node pairs that are separated by the clustering.

\subsubsection{A simple example}

For illustration purpose, inspired by the pedagogical approach of \cite{Waller2017}, we consider the following simple system with $\textbf{x}\in \mathbb{R}^{n_{\textbf{x}}=9}$ and $\textbf{y}\in \mathbb{R}^{n_{\textbf{y}}=4}$:

\begin{align}
   \textbf{H} = 0.25 \times \begin{bmatrix} 
1 & 1 & 1 & 1 & 0 & 0 & 0 & 0 & 0   \\
0 & 1 & 1 & 1 & 0 & 1 & 0 & 0 & 0 \\
0 &  0 & 0 & 1 & 0 & 1 & 1 & 1 & 0  \\
0 & 0 & 0 & 0 & 1 & 1 & 0 & 1 & 1  \label{eq:H_dance}
\end{bmatrix},
\end{align}
where the magnitude of non-zero $\textbf{H}$ entries is assumed constant for simplicity, which leads to the associated state-observation transformation function:
\begin{equation}\label{eq: xy}
\begin{array}{rcl}
    0.25(\textcolor{red}{\textbf{x}_0}+\textcolor{red}{\textbf{x}_1}+\textcolor{red}{\textbf{x}_2}+\textcolor{red}{\textbf{x}_3}) &=& \textbf{y}_0\\
    0.25(\textcolor{red}{\textbf{x}_1}+\textcolor{red}{\textbf{x}_2}+\textcolor{red}{\textbf{x}_3}+\textcolor{blue}{\textbf{x}_5}) &=& \textbf{y}_1\\
    0.25(\textcolor{red}{\textbf{x}_3}+\textcolor{blue}{\textbf{x}_5}+\textcolor{blue}{\textbf{x}_6}+\textcolor{blue}{\textbf{x}_7}) &=& \textbf{y}_2\\
    0.25(\textcolor{blue}{\textbf{x}_4}+\textcolor{blue}{\textbf{x}_5}+\textcolor{blue}{\textbf{x}_7}+\textcolor{blue}{\textbf{x}_8}) &=& \textbf{y}_3. 
    \end{array}
\end{equation}

\begin{figure}
  \centering
        \includegraphics[width=2.6in]{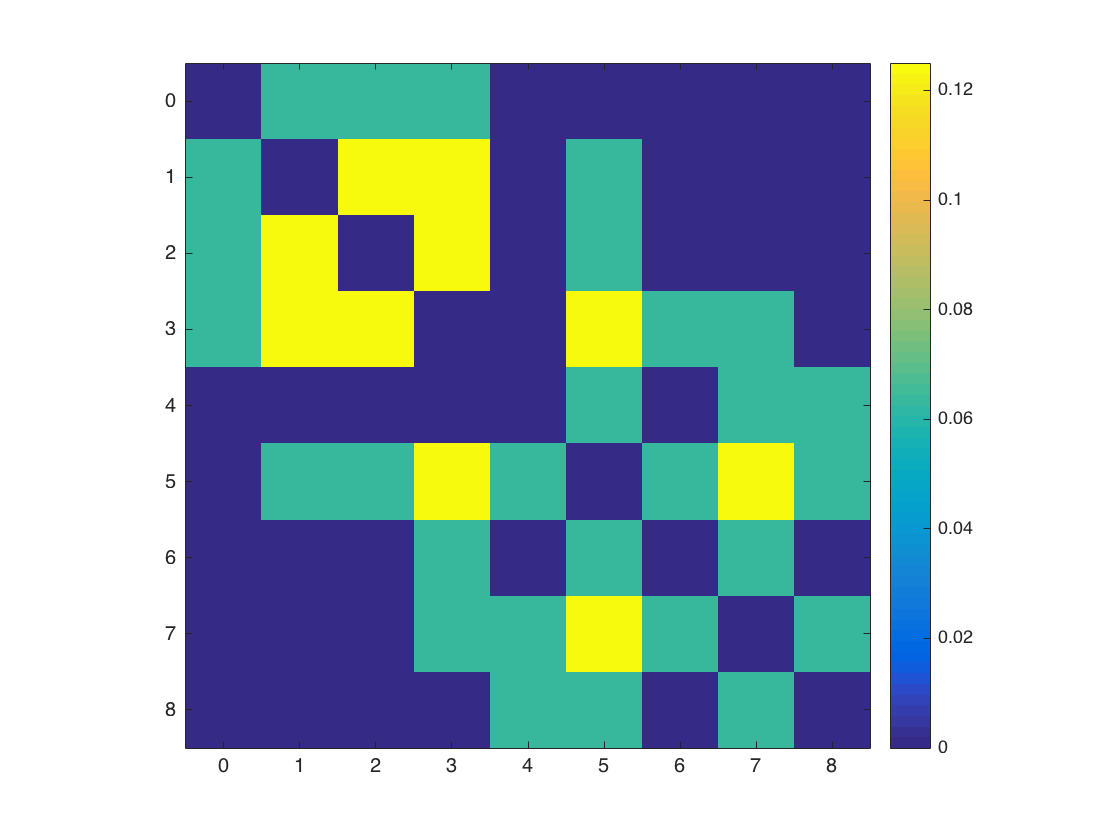} 
        \includegraphics[width=3.0in]{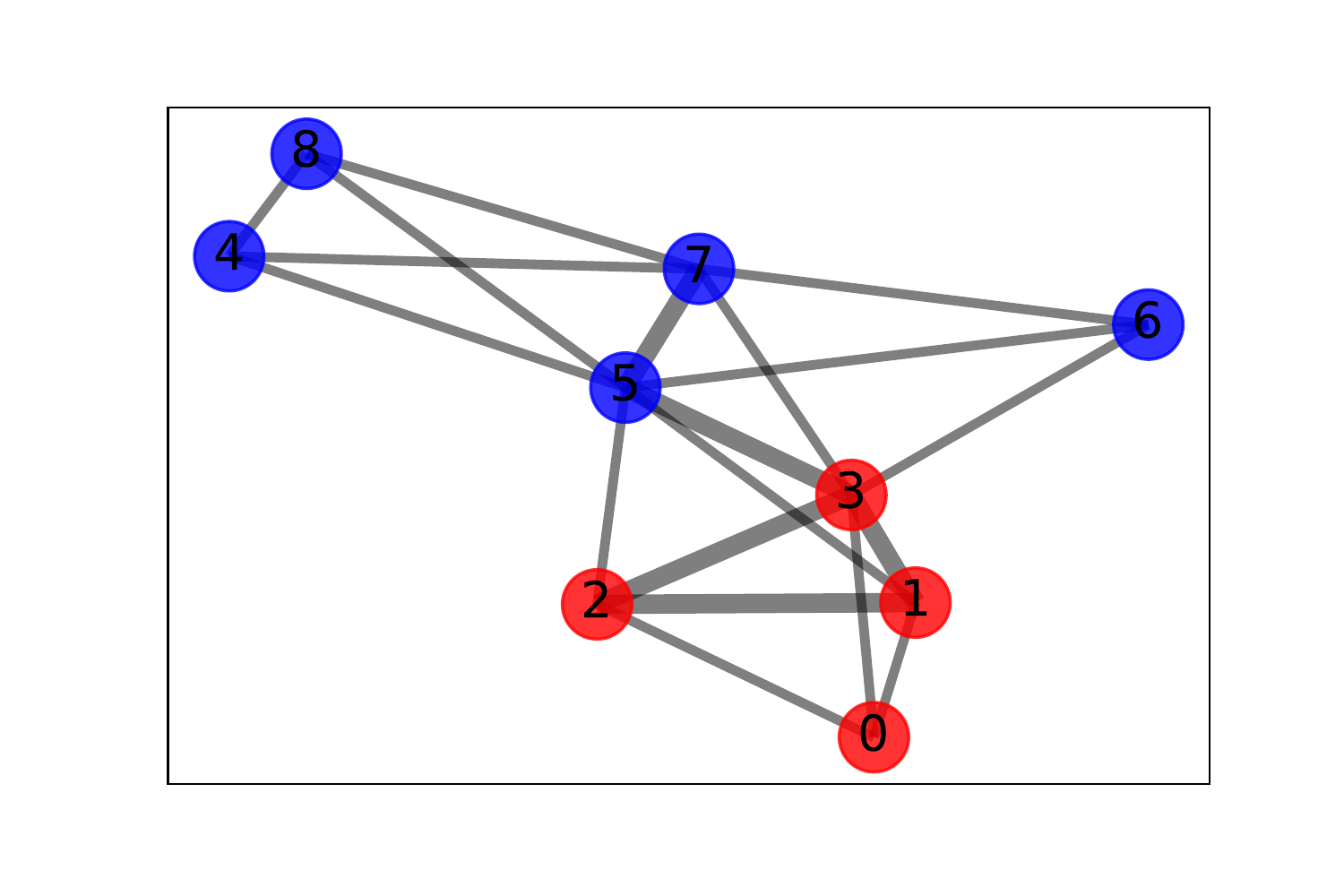} \\
     (a) \hspace{6cm} (b)
     \caption{(a) Observation-based state adjacency matrix obtained from the transformation operator $\textbf{H}$ in Eq. (\ref{eq:H_dance}). (b) Corresponding  network identified by the community detection algorithm. The graph edge weights (measure of the strength of observation-based state connections) are represented by their widths. }
  \label{fig:gaph_ex}
\end{figure}

The obtained observation-based adjacency matrix is represented in Fig. \ref{fig:gaph_ex}(a) and is quite sparse with only $m=19$ edges.
The clustering result obtained by the Fluid community detection algorithm is illustrated in Fig. \ref{fig:gaph_ex}(b). Two communities (red and blue colors) of state variables could be identified, where points in each community are tightly connected. In particular some intra-cluster nodes with strong connections (eg.  $\{\textcolor{red}{\textbf{x}_1},\textcolor{red}{\textbf{x}_2}\}$ or $\{\textcolor{blue}{\textbf{x}_5},\textcolor{blue}{\textbf{x}_7}\}$) are well identified by the algorithm, in accordance with the large values of the adjacency matrix. However,  connections across clusters can also be found, for example the connection $\{\textcolor{red}{\textbf{x}_3},\textcolor{blue}{\textbf{x}_5}\}$. These inter-cluster connections, are still managed by the algorithm. In fact, an output partition of perfect (noise free) subsets can hardly be obtained in real application problems.

After the identification of the state clusters, we need to associate each one with an ensemble of observations. 
As discussed previously, difficulty appears for observations with domains of dependence spanning across multiple clusters.
In this case, it is necessary to operate some data preprocessing. For instance, the assignment of $\{\textbf{y}_0\}$ and $\{\textbf{y}_3\}$ respectively to the first (red) and second (blue) state community is without ambiguity while both $\{\textbf{y}_1\}$ and $\{\textbf{y}_2\}$ are overlapped by the two communities. Dealing with this type of overlapping in the observation partition is therefore crucial for the covariance tuning.

\subsection{Dealing with inter-cluster observation region of dependence for assimilation}

Assuming that a $p$-cluster structure $\mathcal{C}=\{\textbf{x}^1,\textbf{x}^2,\ldots, \textbf{x}^p \}$ {is provided by the community detection algorithm}, we should assign, for each cluster $\textbf{x}^i$, an associated observation subset $\textbf{y}^i$, in order to perform local covariance tuning later on. As we will see in the following, while the partition $\mathcal{C}=\{\textbf{x}^1,\textbf{x}^2,\ldots, \textbf{x}^p \}$ of $\textbf{x}$ will remain the same, the partition of the observations $\{\textbf{y}^1,\textbf{y}^2,\ldots, \textbf{y}^p \}$ will be constructed on a subvector of observations ${\tilde{\textbf{y}}}\in \mathbb{R}^{n_{\tilde{\textbf{y}}}\leq n_{\textbf{y}}}$ or on a modified vector of observations $\hat{\textbf{y}} \in \mathbb{R}^{n_{\textbf{y}}}$. 
In this work, we propose two alternative methods, named ``observation reduction'' and ``observation adjustment'', providing appropriate observation subsets associated with each state cluster. 

\subsubsection{Observation reduction}
  Applying this strategy, the observation components $\textbf{y}_k$ with connections to several state variable clusters must be identified and cancelled, i.e. all observations such that,
  
  \begin{equation}
      \frac{\partial \textbf{y}_{k=0,\ldots,n_{\textbf{y}}-1}}{\partial \textbf{x}^i_{l=1,\ldots,n_{\textbf{x}}^i,i=1,\ldots,p}} \neq 0,
  \end{equation}
  for more than a single cluster, must be withdrawn from the assimilation procedure.
    
    
    Nevertheless, we emphasize that these observation data can still be used later on for evaluating the posterior estimation $\textbf{x}_a$ in the data assimilation procedure. Back to Eq. (\ref{eq: xy}), the  observations $\{\textbf{y}_1\}$ and $\{\textbf{y}_2\}$ are voluntarily excluded to perform the covariance correction, i.e. the tuning will be performed with only two clusters of subvectors: 
\begin{equation}\label{eq: Y1}
\begin{array}{rcl}
        \textbf{x}^1 &=& \{\textbf{x}_{k= 0,\ldots, 3} \}, \quad \tilde{\textbf{y}}^1 = \{\textbf{y}_0\}, \\
        \textbf{x}^2 &=& \{\textbf{x}_{k = 4, \ldots, 8} \},\quad \tilde{\textbf{y}}^2 = \{\textbf{y}_3\}.
    \end{array}
    \end{equation}
    The reduced global state-observation operator $\tilde{\textbf{H}}$ thus becomes
    
\begin{equation}\label{eq:H_reduction} 
 \tilde{\textbf{H}} =  
 0.25 \times \begin{bmatrix} 
1 & 1 & 1 & 1 & 0 & 0 & 0 & 0 & 0    \\
0 & 0 & 0 & 0 & 1 & 1 & 0 & 1 & 1 
\end{bmatrix}.
\end{equation}

\subsubsection{ Observation adjustment}    
 Here, the idea is to modify the observation data dependent on multiple clusters, in order to  simply keep its strongest dependence to a single cluster.   
    This way, each observation will be assigned to only one subset of state variables based on the state-observation mapping. This is done by substracting from the original observation value, the contribution of the surplus quantity related to the other clusters. We rely on the values of the background states to evaluate those surpluses. If more than one background state sample is available (that will be the case  in the next section), the expected value of the background ensemble is used instead.
    
    For example, if $\{\textbf{y}_l \}$ has stronger ties to $\textbf{x}^j$, then it should be readjusted as:
    \begin{align}
        \hat{\textbf{y}}_l = \textbf{y}_l - \sum_{i=1,\ldots,p,i \neq j}\sum_{k \,|\, \textbf{x}_{k} \in \textbf{x}^{i} } \textbf{H}_{l,k} \mathbb{E}_b[\textbf{x}_{k}],
    \end{align}
where $\mathbb{E}_b[.]$ denotes the {\em empirical} expected value based on the prior background ensemble at hand. This approach leads to an adjusted Jacobian matrix $\hat{\textbf{H}}$ that induces adjacency matrix with no overlapped domains. This is obviously an approximation due to the averaged operator. In fact, there are two error sources, a main one coming from the prior backgroup measure and another one due to the sampling error. We will see examples 
in section \ref{sec:Experimental results}
    
    Applied to the example, $\{\textbf{y}_1\}$ and $\{\textbf{y}_2\}$ can be respectively adjusted to belong to the first and the second cluster. With the help of background state $\textbf{x}_b$, Eq. (\ref{eq: xy}) can be adjusted to:
    \begin{equation}\label{eq: xy_reshape}
        \begin{array}{l}
    0.25(\textcolor{red}{\textbf{x}_0}+\textcolor{red}{\textbf{x}_1}+\textcolor{red}{\textbf{x}_2}+\textcolor{red}{\textbf{x}_3}) = \hat{\textbf{y}}_0 = \textbf{y}_0\\
    0.25(\textcolor{red}{\textbf{x}_1}+\textcolor{red}{\textbf{x}_2}+\textcolor{red}{\textbf{x}_3}) = \hat{\textbf{y}}_1 =\textbf{y}_1 - 0.25 \textcolor{blue}{\mathbb{E}_b[\textbf{x}_{5}]}\\
    0.25(\textcolor{blue}{\textbf{x}_5}+\textcolor{blue}{\textbf{x}_6}+\textcolor{blue}{\textbf{x}_7}) =  \hat{\textbf{y}}_2 = \textbf{y}_2 -0.25 \,\textcolor{red}{\mathbb{E}_b[\textbf{x}_{3}]}\\
    0.25(\textcolor{blue}{\textbf{x}_4}+\textcolor{blue}{\textbf{x}_5}+\textcolor{blue}{\textbf{x}_7}+\textcolor{blue}{\textbf{x}_8}) = \hat{\textbf{y}}_3 = \textbf{y}_3. 
        \end{array}
    \end{equation}
Thus the new operator 
can be written as:

\begin{align}
   \hat{\textbf{H}}_{\textrm{adjustment}} 
   & = 0.25 \times  \begin{bmatrix} 
1 & 1 & 1 & 1 & 0 & 0 & 0 & 0 & 0   \\
0 & 1 & 1 & 1 & 0 & 0 & 0 & 0 & 0 \\
0 &  0 & 0 & 0 & 0 & 1 & 1 & 1 & 0  \\
0 & 0 & 0 & 0 & 1 & 1 & 0 & 1 & 1 \label{eq:H_reshape}
\end{bmatrix}.
\end{align}
For real applications, one may envision a mixture of these two approaches.

\section{Desroziers \& Ivanov covariance tuning algorithm }
\label{sec:Desroziers Ivanov tuning algorithm}

Now that we have localized our system based on the state-observation linearized measure, and thanks to graph clustering methods, we are going to explain how we proceed with our local covariance tuning within each subdomain.

\subsection{Error covariance diagnosis and tuning}

The \cite{Desroziers01} tuning algorithm (\di)  was first proposed and applied in the meteorological science at the beginning of the 21st century. This method is based on the diagnosis and verification of innovation quantities  and has been widely applied in meteorology and geoscience (e.g. \cite{Hoffman2013}). Consecutive works have been carried out to improve its performance and feasibility in problems of large dimension such as the study of \cite{Chapnik2006}.  

It was proven in \cite{Talagrand99} and \cite{Desroziers01} that, under the assumption of perfect knowledge of the covariance matrices $\textbf{B}$ and $\textbf{R}$, 
the following equalities are perfectly satisfied in a  3D-VAR assimilation system:
\begin{align}
	\mathbb{E}\left [J_b(\textbf{x}_a) \right ] & = \frac{1}{2} \mathbb{E}\left [(\textbf{x}_a-\textbf{x}_b)^T\textbf{B}^{-1}(\textbf{x}_a-\textbf{x}_b) \right ]\\
	& =\frac{1}{2} \textrm{Tr}(\textbf{K}\textbf{H}), \label{eq:Jb} \notag
\end{align} 
\begin{align}
	\mathbb{E} \left [J_o(\textbf{x}_a)  \right ] & = \frac{1}{2} \mathbb{E}\left [(\textbf{y}-\textbf{H}\textbf{x}_b)^T\textbf{R}^{-1}(\textbf{y}-\textbf{H}\textbf{x}_b) \right ]\\
	& =\frac{1}{2} \textrm{Tr}(\mathcal{I}-\textbf{H}\textbf{K}), \label{eq:Jo} \notag
\end{align}
where $ \textbf{x}_a $ is the output of a 3D-VAR algorithm with a linear observation operator $\textbf{H}$. \smallskip 
In practice, this is seldomly the case and one has to deal with imperfect knowledge of the covariance matrices. Nonetheless, if we assume that the correlation structures of these matrices are well known, then it is possible to iteratively correct their magnitudes.
Using the two indicators
\begin{equation}
	s_{b,q}=\frac{2J_b(\textbf{x}_a)}{\textrm{Tr}(\textbf{K}_q \textbf{H})},
\end{equation}
\begin{equation}
	s_{o,q}=\frac{2J_o(\textbf{x}_a)}{\textrm{Tr}(\mathcal{I}-\textbf{H}\textbf{K}_q)}, \label{eq:DI01}
\end{equation}
where $q$ is the current iteration, 
the objective of the \di ~tuning method is to adjust the ratio between the weighting of $\textbf{B}^{-1}$ and $\textbf{R}^{-1}$ without modifying their correlation structure:
\begin{equation}
	\textbf{B}_{q+1}=s_{b,q} \textbf{B}_q, \quad \textbf{R}_{q+1}=s_{o,q} \textbf{R}_q.
\end{equation}
These two indicators act as scaling coefficients, modifying the error variance magnitude. We remind that both the reconstructed state $\textbf{x}_a$ and the gain matrix $\textbf{K}_q$ depend on $\textbf{B}_q$, $\textbf{R}_q$ and thus on the iterative coefficients $s_{b,q}$, $s_{o,q}$. The application of this method in subspaces where matrices $\textbf{B}$ and $\textbf{R}$ follow block-diagonal structures has also been discussed in \cite{Desroziers01}.

Compared to other posterior diagnosis or iterative methods, e.g. \cite{Desroziers2005}, \cite{cheng2019}, no estimation of full matrices is needed and only the estimation of two scalar values ($J_b,J_o$) is required in \di. Therefore, this method could be more suitable when the available data is limited. Another advantage relates to the computational cost of this method as DI01 requires only the computation of matrices trace which can be evaluated in efficient ways. 

According to \cite{Chapnik2006}, the convergence of $s_b$ and $s_o$ can be very fast, especially in the ideal case where the correlation patterns of $\textbf{B}$ and $\textbf{R}$ are perfectly known. Under this assumption, \cite{Chapnik2006} proved DI01 is equivalent to a maximum-likelihood parameter tuning. In addition, large iteration number is not required as the first iteration could already provide a reasonably good estimation of the final result.

\subsection{Application of \di ~in subspaces}
The application of data assimilation algorithms, as well as the full  observation matrix diagnosis has been discussed in \cite{Waller2017}.   Following the notation of their paper, we  introduce the binary selection matrix $\Phi_{x}^i, \Phi_{y}^i$ of the $i^{th}$ subvector with
\begin{align}
    \textbf{x}^{i}=\Phi_{x}^i \textbf{x}, \quad \textbf{y}^{i}=\Phi_{y}^i \textbf{y}
\end{align}
where $i$ is the index of the subspace. The data assimilation in the subspace, as well as localized covariance tuning could be easily expressed using the standard formulation with projection operators  $\Phi_{x}^i$ and $ \Phi_{y}^i$.

Given the example of the first pair of state and observation subsets in the case of Fig. \ref{fig:gaph_ex}, we have
\begin{align}
    \Phi_{x}^1 =  \begin{bmatrix} 
1 & 0 & 0 & 0 & 0 & 0 & 0 & 0 & 0   \\
0 & 1 & 0 & 0 & 0 & 0 & 0 & 0 & 0 \\
0 &  0 & 1 & 0 & 0 & 0 & 0 & 0 & 0  \\
0 & 0 & 0 & 1 & 0 & 0 & 0 & 0 & 0 \label{eq:x_H_reshape}
\end{bmatrix}.
\end{align}
In the case of data reduction strategy (Eq.\ref{eq: Y1}),
\begin{align}
    \Phi_{y,\textrm{reduction}}^1 =  \begin{bmatrix} 
1 & 0 & 0 & 0    \label{eq:phi_reduction}
\end{bmatrix},
\end{align}
while for data adjustment strategy, 
\begin{align}
    \Phi_{y,\textrm{adjustment}}^1 =  \begin{bmatrix} 
1 & 0 & 0 & 0    \\
0 & 1 & 0 & 0   \label{eq:phi_reshape}
\end{bmatrix}.
\end{align}

The error covariances matrix $\textbf{B}^{i} $ (resp. $\textbf{R}^{i} $) associated to $\textbf{x}_b^{i}$ (resp. $\textbf{y}^{i}$) can be written as:
\begin{align}
   \textbf{B}^{i} =   \Phi_{x}^i \textbf{B} \Phi_{x}^{i,T}, \quad \textbf{R}^{i} =   \Phi_{y}^i \textbf{R} \Phi_{y}^{i,T}.
\end{align}
Therefore, the associated analyzed subvector $\textbf{x}_a^{i}$ could be obtained by applying data assimilation procedure using $\Big( {\textbf{x}_b}^{i}, \textbf{y}^{i}, \textbf{B}^{i}, \textbf{R}^{i} \Big)$. We remind that, due to the cross-community noises (i.e. the updating of ${\textbf{x}_b}^{i} $ may not only depend on ${\textbf{y}_b}^{i} $ in the global data assimilation system), we don't necessarily have 
\begin{align}
    {\textbf{x}}_a^{i}=\Phi_x^i \textbf{x}_a.
\end{align}
For more details of decomposition formulations, the interested readers are referred to \cite{Waller2017}.

\subsection{Localized \di ~tuning }
Our objective for implementing localized covariance tuning algorithms is to gain a finer diagnosis and correction on the covariance computation. The local DI01 diagnosis in $\Big( {\textbf{x}_b}^{i}, {\textbf{y}}^{i} \Big) $ can be expressed as:
\begin{align}
	\mathbb{E}\left [ J_b(\textbf{x}_a^i)\right ] & = \mathbb{E}\left [(\textbf{x}_a^i-\textbf{x}_b^i)^T({\textbf{B}}^{i})^{-1}(\textbf{x}_a^i-\textbf{x}_b^i) \right ] \notag \\
	& =\frac{1}{2} \textrm{Tr}({\textbf{K}}^{i} {\textbf{H}}^{i}) , \\
	\mathbb{E}\left [J_o(\textbf{x}_a^i)\right ]&=\mathbb{E}\left [({\textbf{y}}^{i}-{\textbf{H}}^{i} \textbf{x}_b^i)^T ({\textbf{R}}^{i})^{-1}({\textbf{y}}^{i}-{\textbf{H}}^{i} \textbf{x}_b^i)\right ] \notag \\
	& =\frac{1}{2} \textrm{Tr}(\mathcal{I}^i-{\textbf{H}}^{i} {\textbf{K}}^{i}),
\end{align}
where the optimization functions $ J_b$ and $J_o$, as well as the localized gain matrix $ {\textbf{K}}^{i}$ have also been adjusted in these subspaces.  
The identity matrix $\mathcal{I}^i$ is of the same dimension as ${\textbf{B}}^{i}$. 

We can also define the local tuning algorithm:
\begin{align}
	s_{b,q}^{i}&=\frac{2J_b(\textbf{x}_a^i)}{\textrm{Tr}({\textbf{K}_q}^{i}{\textbf{H}}^{i})} \label{eq:local_di_b},\\
	s_{o,q}^{i}&=\frac{2J_o(\textbf{x}_a^i)}{\textrm{Tr}(\mathcal{I}^i-{\textbf{H}}^{i}{\textbf{K}_q}^{i})}, \label{eq:local_di_o}\\
	{\textbf{B}}^{i}_{q+1}&=s_{b,q}^{i} {\textbf{B}}^{i}_{q}, \\
	{\textbf{R}}^{i}_{q+1}&=s_{o,q}^{i} {\textbf{R}}^{i}_{q}. 
\end{align}
{The iterative process is repeated $q_\textrm{max}^i$ times, based on some \textit{a priori} maximum number of iterations or some stopping criteria monitoring the rate of change. }  The approach provides a local correction within each cluster thanks to a  multiplicative coefficient. 
This way, the covariance tuning is more flexible than a global approach relying on two coefficients $(s_b, s_o)$ only.

However, if the updating is performed in each subspace (i.e. correction only on the sub-matrices $ {\textbf{B}}^{i}, {\textbf{R}}^{i}$), then the adjusted $\textbf{B}$ and $\textbf{R}$ are not guaranteed to be positive-definite. In order to circumvent this problem, we keep the correlation structure of $(\textbf{C}_{\textbf{B}}$ and $\textbf{C}_{\textbf{R}})$ fixed. We remind that a covariance matrix $\textbf{Cov}$ (of random vector $\textbf{x}$), which is by its nature positive semi-definite, can be decomposed into its variance and correlation structures as:
\begin{equation*}
      \textbf{Cov} = \textbf{D}^{1/2} \textbf{C} \hspace{1mm} \textbf{D}^{1/2},  
\end{equation*}
where $\textbf{D}$ is a diagonal matrix of the state error variances,
and $\textbf{C}$ is the correlation matrix. By correcting the variance in each subspace only through the diagonal matrices
$(\textbf{D}_{\textbf{B}}^{i},\textbf{D}_{\textbf{R}}^{i})$, the positive definiteness of $\textbf{B}$ and $\textbf{R}$ is thus guaranteed as the correlation structure remains invariant, cf. Algorithm \ref{algo:1}.
\vspace{5mm}

\begin{algorithm}[]
\caption{Localization and updating of $\textbf{B}$ and $\textbf{R}$ with cluster-based implementation of DI01 algorithm. }
Inputs: \\
Background state: $\textbf{x}_b$ \\
Observation data: $\textbf{y}$\\
Initially guessed matrix: $\textbf{B}, \textbf{R}$\\
Jacobian matrix: $\textbf{H}$\\
\vspace{2mm}
Algorithm:
Community detection using $\textbf{H}$ with given or detected community number~$p$\\
\For(){$i$ from $1$ to $p$:}
 {
    {Loading of subvectors ${\textbf{x}_b}^{i}, {\textbf{y}}^{i} $ and associated covariance matrices ${\textbf{B}}^{i},{\textbf{R}}^{i}$. } \smallskip \\
%
Initializing covariance matrices:
\begin{align}
\textbf{B} &\leftarrow \textbf{B}_{\textbf{A}}= \sigma_{b,\textbf{A}}^2\, \textbf{C}_{\textbf{B}} \notag \\
\textbf{R} &\leftarrow \textbf{R}_{\textbf{A}}= \sigma_{o,\textbf{A}}^2\, \textbf{C}_{\textbf{R}} \notag
\end{align}

    \For(){$q$ from $1$ to $q_{\textrm{max}}$}
    {
    calculation of $\left \{ s_{b,q}^{i},s_{o,q}^{i}\right \}$, Eq.~(\ref{eq:local_di_b}-\ref{eq:local_di_o})
    }
Updating of full covariance matrices from blockwise tuned covariance in current cluster:
    \begin{align}
\textbf{B} &\leftarrow ({\textbf{D}}_\textbf{B}^{i})^{1/2} \textbf{B} \hspace{1mm} ({\textbf{D}}_\textbf{B}^{i})^{1/2} \notag \\
\textbf{R} &\leftarrow ({\textbf{D}}_\textbf{R}^{i})^{1/2} \textbf{R} \hspace{1mm} ({\textbf{D}}_\textbf{R}^{i})^{1/2} \notag
\end{align}
where ${\textbf{D}}_\textbf{B}^{i}$ and ${\textbf{D}}_\textbf{R}^{i}$ are diagonal matrices defined as:
\begin{align}
({\textbf{D}}_\textbf{B}^{i})_{j,j} =\left\{
\begin{array}{c l}	
     & \prod_{q=1}^{q_\textrm{max}} s_{b,q}^{i} \quad \text{if}\; \{\textbf{x}_j\} \subset \textbf{x}^i \\
     & 1 \quad \textrm{otherwise}
\end{array}\right. \notag \\
({\textbf{D}}_\textbf{R}^{i})_{l,l}=\left\{
\begin{array}{c l}	
     & \prod_{q=1}^{q_\textrm{max}} s_{o,q}^{i} \quad \textrm{if}\; \{\textbf{y}_l\} \subset \textbf{y}^i\\
     & 1 \quad \textrm{otherwise}. \notag
\end{array}\right.
\end{align}
 }

 \vspace{2mm}
outputs: Improved error covariances
\label{algo:1}
\end{algorithm}

\section{Illustration with numerical experiments }
\label{sec:Illustration with numerical experiments}
\subsection{Test case description}
\label{sec:Test case description}
\subsubsection{Construction of $\textbf{H}$}
\label{sec:Construction of H}

Here we describe our numerical experiments which combine -- a real algorithm of community detection, -- an artificially created state-observation mapping and -- the implementation of covariance tuning method in the underlying subspaces. The numerical tests are carried out in the same spirit as in the paper by  \cite{Waller2017} but for a larger system. A sparse Jacobian matrix $\textbf{H}$ reflecting the clustering of the state-observation mapping is generated; the components of which are then randomly mixed in order to hide any particular structural pattern. The dimension of the state space is set to be 100, $\textbf{x} \in \mathbb{R}^{n_{\textbf{x}}=100}$, while the dimension of the observation space is set to be $\textbf{y} \in \mathbb{R}^{n_{\textbf{y}}=50}$.  We consider a case for which the state-observation mapping $\textbf{H}$ reflects community structures.  For this reason, we construct \textit{a priori} two (this choice is arbitrary) subsets of observations each relating mainly to only one subset of state variables. In order to be as general as possible,
we consider $\left | \textbf{x}^1 \right | = \left | \textbf{x}^2 \right |=50$ and $\left | \textbf{y}^1 \right | = \left | \textbf{y}^2 \right |=25$.
For the sake of simplicity, the observation operator $\textbf{H}$ (of dimension $[50 \times 100]$) is randomly filled with binary elements, forming a dominant blockwise structure with some extra-block non-zero terms. The latter is done in order to mimic realistic problems, i.e. some perturbations are introduced in the form of cross-communities perturbations, therefore the two communities are not perfectly separable.

The background/observation vectors and Jacobian matrix are then randomly shuffled in a coherent manner in order to hide the cluster structure to the community detection algorithm, as for the adjacency matrix in Fig. \ref{fig:ad_matrices}(a).
More specifically, the state-observation mapping is constructed as follows:
we use a binomial distribution with two levels of success probability:
\begin{align}
Pr(\textbf{H}_{i,j} & = 1) = \\
& \left\{
    \begin{array}{ll}
        15\% \quad \textrm{if} \quad \textbf{x}_i \in \textbf{x}^{1} \quad \textrm{and} \quad y_j \in \textbf{y}^{1}  \\
        15\% \quad \textrm{if} \quad \textbf{x}_i \in \textbf{x}^{2} \quad \textrm{and} \quad  y_j \in \textbf{y}^{2}  \\
        1\% \quad \textrm{otherwise} \quad \textrm{(perturbations)}. \notag
    \end{array}
\right.
\end{align}

In the following tests, exact and assumed covariance magnitudes will be changed but we will always keep the same choice of Jacobian $\textbf{H}$. The community detection, remaining also invariant for all Monte Carlo tests, is provided by the Fluid community-detection algorithm. As explained previously, there is a particular interest to apply \di ~in the case of limited access to data (i.e. small ensemble size of ($\textbf{x}_b,\textbf{y}$)).

In these twin experiments, the prior errors are assumed to follow the distribution of correlated Gaussian vectors:
\begin{align}
    \epsilon_b &=\textbf{x}_b-\textbf{x}_t \sim \mathcal{N}(0^{n_{\textbf{x}}=100},\textbf{B}_\textrm{E}),\\
    \epsilon_y &= \textbf{y} - \textbf{H}\textbf{x}_t  \sim \mathcal{N}(0^{n_{\textbf{y}}=50},\textbf{R}_\textrm{E}),
\end{align}

where $\textbf{B}_E, \textbf{R}_E $ denote the chosen exact prior error covariances, {\em hidden from the tuning algorithm}. 
We remind that under the assumption of state independent error and linearity of $\textbf{H}$, the posterior assimilation error, as well as the posterior correction of $\textbf{B}$ and $\textbf{R}$ via \di  ~(regardless of the strategy chosen, i.e. data reduction or data adjustment), is independent of the theoretical value of $\textbf{x}_t$ but only depends on prior errors (i.e. $\textbf{x}_t-\textbf{x}_b$ and $\textbf{y}-\textbf{H}\textbf{x}_t$).

\subsubsection{Twin experiments setup}
\label{sec:Twin experiments set up}
In order to reflect the construction of $\textbf{H}$, we suppose that the {\em exact} error deviation, {\em hidden from the tuning algorithm} (respectively denoted by $\sigma_{b,E}^i,\sigma_{o,E}^i$) are constant in each cluster, so for instance we have:
\begin{align}
    \textrm{if} \quad \{\textbf{x}_u, \textbf{x}_v \} \subset \textbf{x}^i, \quad \textrm{then} \quad \sigma_{b,E}^i(\textbf{x}_u) = \sigma_{b,E}^i(\textbf{x}_v). \nonumber
\end{align}
For this numerical experiment, a quite challenging case is chosen with:
\begin{eqnarray}
\sigma_{b,E}^{i=1}(\textbf{x}_u) & = & \sigma_{b,E}^{i=2}(\textbf{x}_v) \nonumber \\
\sigma_{o,E}^{i=1}(\textbf{y}_u) & = & \text{ratio}\times\sigma_{o,E}^{i=2}(\textbf{y}_v), \nonumber
\end{eqnarray}
so that the background error is homogeneous while the observation error is different in the two communities with a fixed ratio (in the following, we will choose $ratio=10$).
However, the correlation structures of the covariance matrices are supposed to be known {\em a priori}, and are assumed to follow a Balgovind structure:
\begin{align}
   (\textbf{C}_{\textbf{B}})_{i,j}=(\textbf{C}_{\textbf{R}})_{i,j} = \left ( 1+\frac{r}{L}\right ) \exp^{-\frac{r}{L}},
\end{align}
 where $r\equiv r(\textbf{x}_i, \textbf{x}_j) = r(\textbf{y}_i, \textbf{y}_j) = |i-j|$ is a pseudo spatial distance between two state variables, 
and the correlation scale is fixed ($L=10$) in the following experiments.
The Balgovind structure is also known as the $\nu=3/2$ Matern kernel, often used in prior error covariance computation in data assimilation (see for example \cite{Poncot2013}, \cite{dance2013}).
%

We remind that the output of all DI01 based approaches depend on the available background and observation data set. We compare three different methods described previously in this paper, differentiated by the notation used for their output covariances:
\begin{itemize}
    \item $(\textbf{B}, \textbf{R})$: implementation of \di~in full space,
    \item $(\tilde{\textbf{B}}, \tilde{\textbf{R}})$: implementation of \di~with graph clustering localization with {\em data reduction} strategy,
   \item $(\hat{\textbf{B}}, \hat{\textbf{R}})$: implementation of \di~with graph clustering localization with {\em data adjustment} strategy.
\end{itemize}

The performance of the covariance tuning with localization is evaluated with a simple scalar criteria involving the Frobenius norm, relative to the standard approach. This indicator/gain may be expressed for the background covariance tuning with the reduction strategy as: 
\begin{align}
    \gamma_{\tilde{\textbf{B}}} = \left ( \Delta_{{\textbf{B}}} -\Delta_{\tilde{\textbf{B}}} \right )/\Delta_{{\textbf{B}}},
\end{align}
with $\Delta_{\cdot}= \mathbb{E}\big[ \|\cdot-\textbf{B}_\textrm{E}\|_F \big]$ representing the expected matrices difference in the Frobenius norm, and similarly for the adjustment strategy and for the observations covariance. 
The larger the gain, the more advantage can be expected from the new approach compared to the standard \di~one.

\begin{figure}
  \centering
     \includegraphics[trim=7.7cm 1cm 4cm 8cm,clip=true,width=3.in]{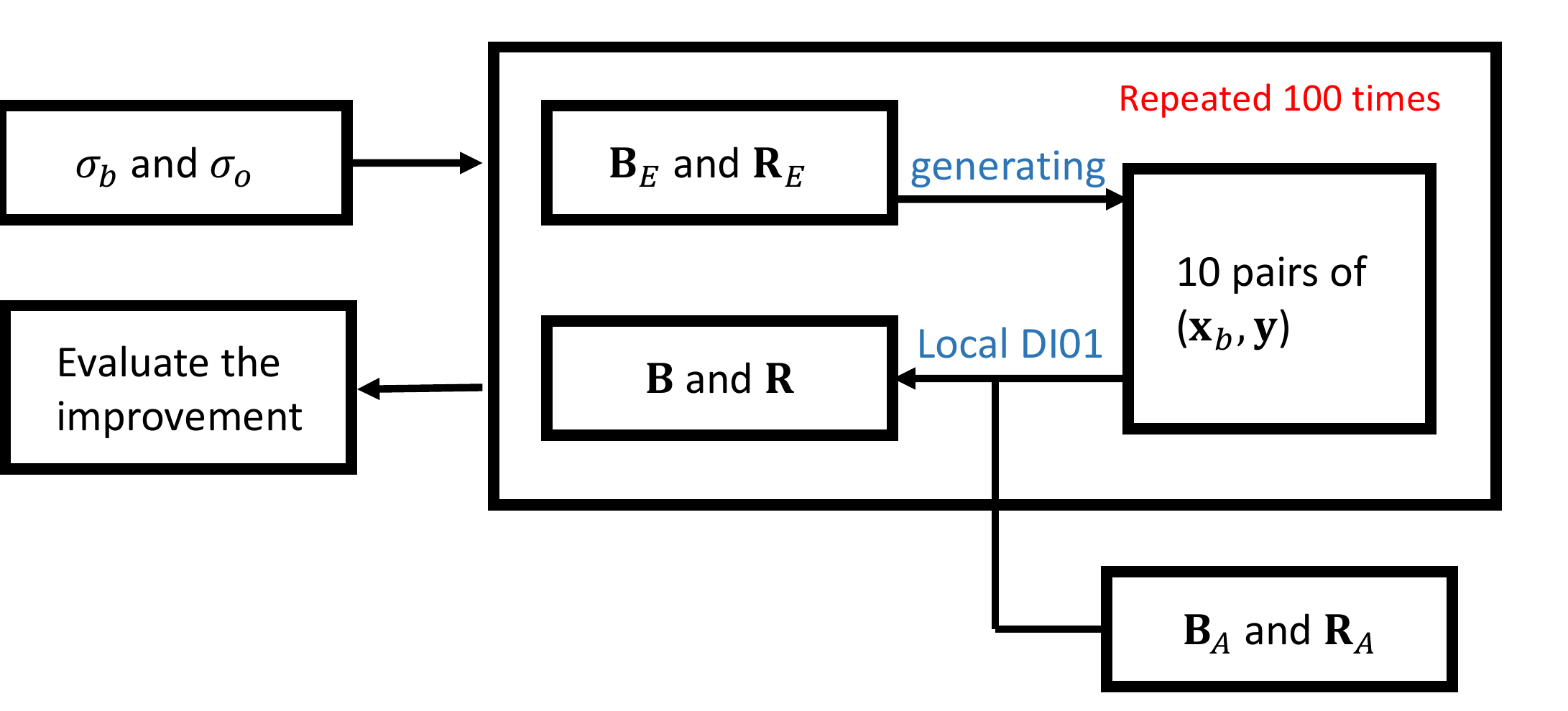}   
     \caption{Flowchart of Monte-Carlo experiments for adaptive data assimilation with fixed parameters: $\sigma_b, \sigma_o, \textbf{B}_A, \textbf{R}_A$.}
  \label{fig:principal_test}
\end{figure}

In the numerical results presented later, empirical expectation of these indicators will be calculated by repeating the tests $100$ times, in a Monte Carlo fashion, for each case of standard deviation parameters. In each Monte Carlo simulation, 10 pairs of background state and observation vector are generated to evaluate the coefficients $s_b^i$ and $s_o^i$ necessary for diagnosing and improving $\textbf{B}$ and $\textbf{R}$ as shown in Fig. \ref{fig:principal_test}.

In order to examine the performance of the proposed approach, we choose to always quantify the assumed prior covariances as;
\begin{align*}
    \textbf{B}_A &= \sigma_{b,A}^2 \times \textbf{C}_{\textbf{B}},  \\
    \textbf{R}_A &= \sigma_{o,A}^2 \times \textbf{C}_{\textbf{R}},
\end{align*}
with $\sigma_{b,A}=\sigma_{o,A}=0.05$.

Meanwhile, the average exact prior error deviation ($\sigma_{b,E}, \sigma_{o,E} = \sqrt{\sigma_{o,1} \sigma_{o,2}}$) varies in the range $([0.025,0.1])$. In other words, we test a range spanning a domain with over-estimation of $100\%$ of error deviation to an under-estimation of $100\%$.  {We remind that the aim of the new approaches is to obtain a more precise estimation of prior covariance structures. }. 

\subsection{Results}
\label{sec:Experimental results}
\begin{figure}
  \centering
     \includegraphics[width=2.in]{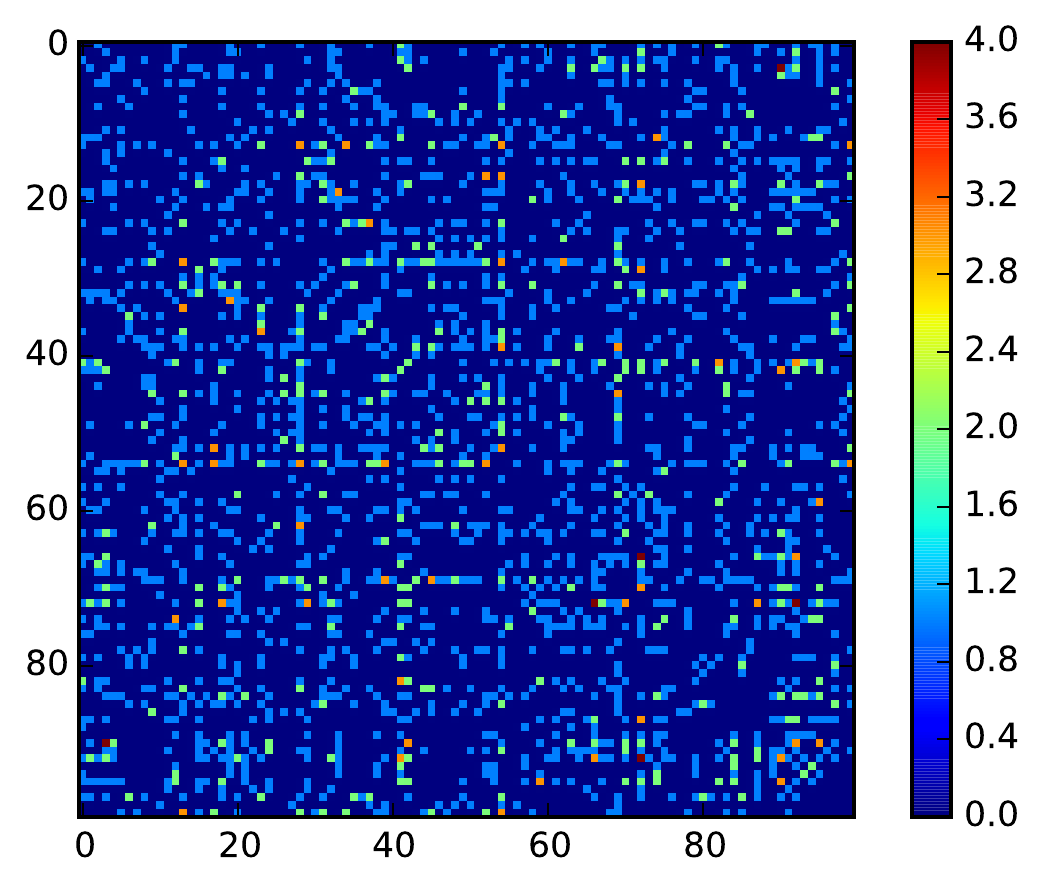} \hspace{1.8cm}
    \includegraphics[width=2.in]{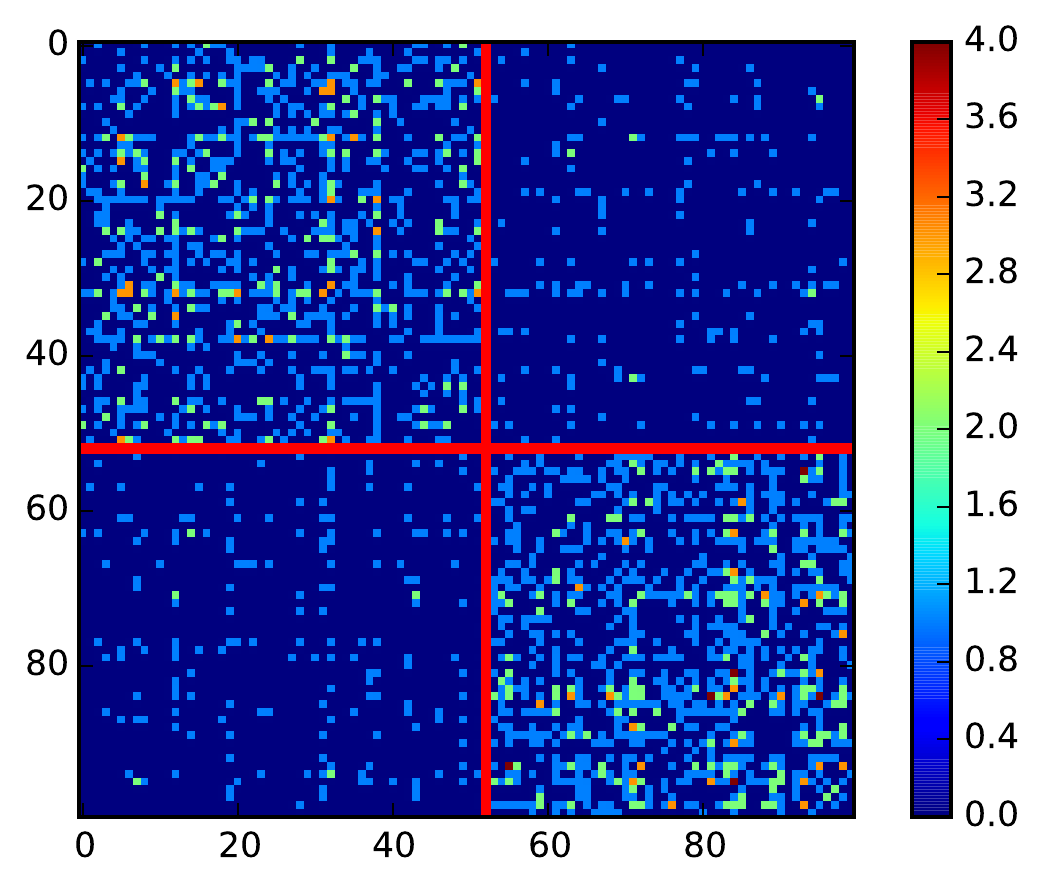} \\
     (a) \hspace{6cm} (b)
     \caption{(a) Original adjacency matrix of a 100 vertex observation-based state network. (b) Vertex ordering by cluster where the 2-cluster structure is evident thanks to the graph clustering algorithm.} 
  \label{fig:ad_matrices}
\end{figure}

\begin{figure}
  \centering
     \includegraphics[width=2.8in]{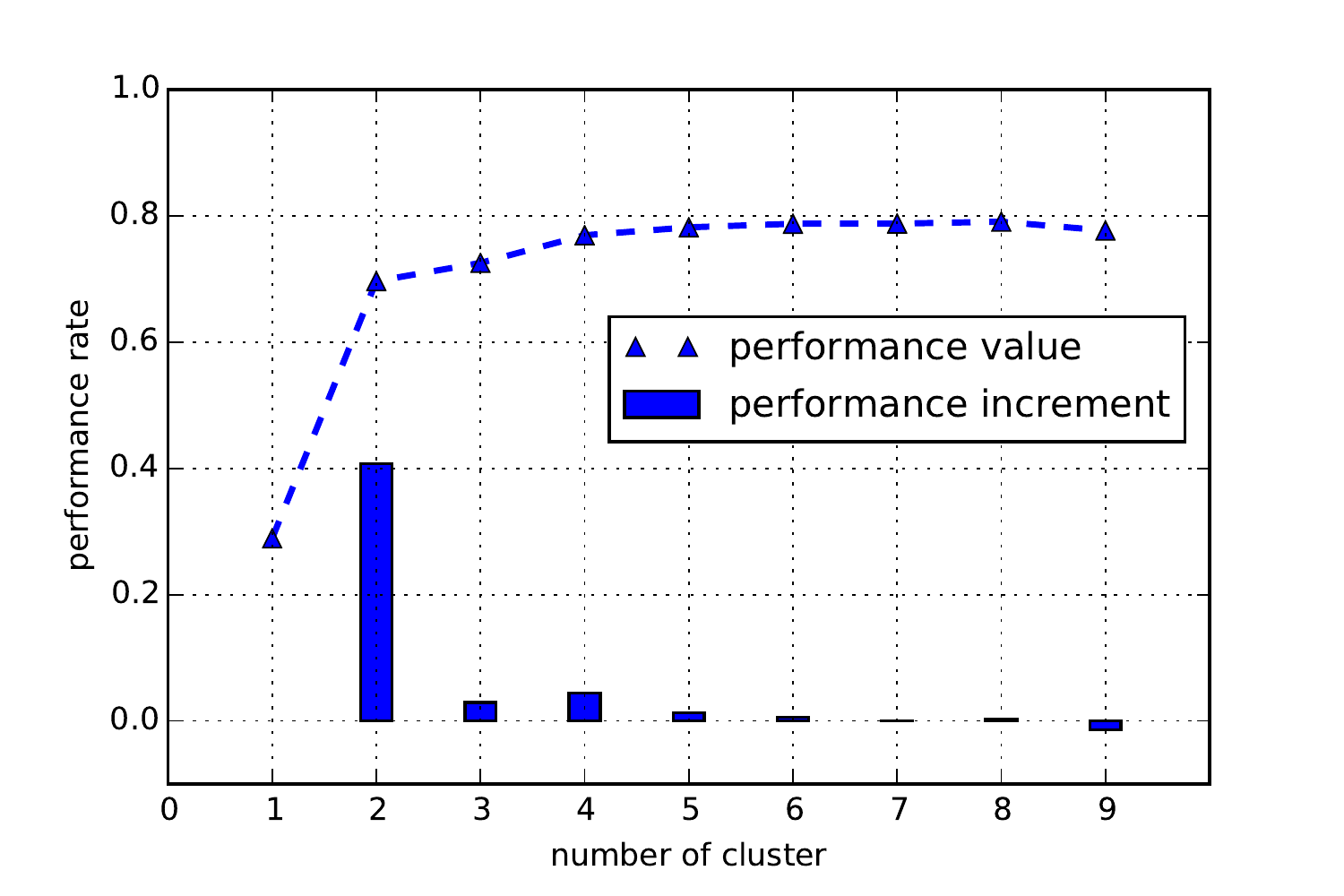} 
     \caption{The evolution of performance value and its increment against the number of communities chosen.}
  \label{fig:coverage}
\end{figure}

{Thanks to the adjacency matrix of the observation-based state network as shown in Fig. \ref{fig:ad_matrices}(a)}, we first apply the Fluid community detection method in the observation-based state network to determine subspaces (communities) in the state space. For real applications, the number of communities is unknown. Here, we apply several times the community detection algorithms with different assumed community number and we evaluate the performance rate (\cite{FORTUNATO2010}) of the obtained partition. It is used as an indicator for finding the optimal community number (as shown in Fig. \ref{fig:coverage}), which is a standard approach for graph problems.

According to the result presented in Fig. \ref{fig:coverage}, we chose to separate the state variables into two subsets, which is the correct number of communities when we simulate the Jacobian matrix $\textbf{H}$. {We emphasize that despite the fact that the $\textbf{H}$ matrix was generated using two clusters, it was not trivial to rediscover them from the observation-based state network, once the information was shuffled and noisy, cf. from (a) to (b) in Fig. \ref{fig:ad_matrices}.} The result of graph partition algorithm of two communities is summarized in Table. \ref{table:1} and Fig. \ref{fig:ad_matrices}(b). 
From Table. \ref{table:1}, we notice that two state variables from the second subset $\textbf{x}^2$ are mistakenly assigned to the first one, $\textbf{x}^1$. The last column of Table. \ref{table:1} shows the total number of observations ($| \textbf{y} |= | \textbf{y}^1 | + |\textbf{y}^2 |$) used in the covariance tuning. Only half of the observations are considered while applying the strategy of data adjustment.

\begin{table*}[t]
  \centering
  \resizebox{0.95\linewidth}{!}{%
\centering
    \begin{tabular}{ | l|cccccc |}
    \hline
     \multirow{2}{*}{\backslashbox{Strategy chosen}{Size of subsets}} &  &    &  Detected & subsets & & \\
    \cline{2-7}  
     & $|\textbf{x}^1|$  & $| \textbf{x}^2 |$  & $ | \textbf{x}|$  & $|\textbf{y}^1|$  & $|\textbf{y}^2|$  & $| \textbf{y} |$\\
    \hline  
    \hline
     {Data reduction} &  52 & 48 & 100 &  12
 & 13 & 25
 \\
     {Data adjustment} &   52 & 48 &  100 &  25
 & 25 & 50 \\
      \hline
  \end{tabular}}
  \caption{Quantification of the community detection algorithm results on the observation-based state network followed by the data reduction and data adjustment strategies. The number of communities (i.e. $k=2$) is set according to the result in Fig. \ref{fig:coverage}. }
  \label{table:1}
\end{table*}

\begin{figure*}
  \centering
     \includegraphics[width=2.65in]{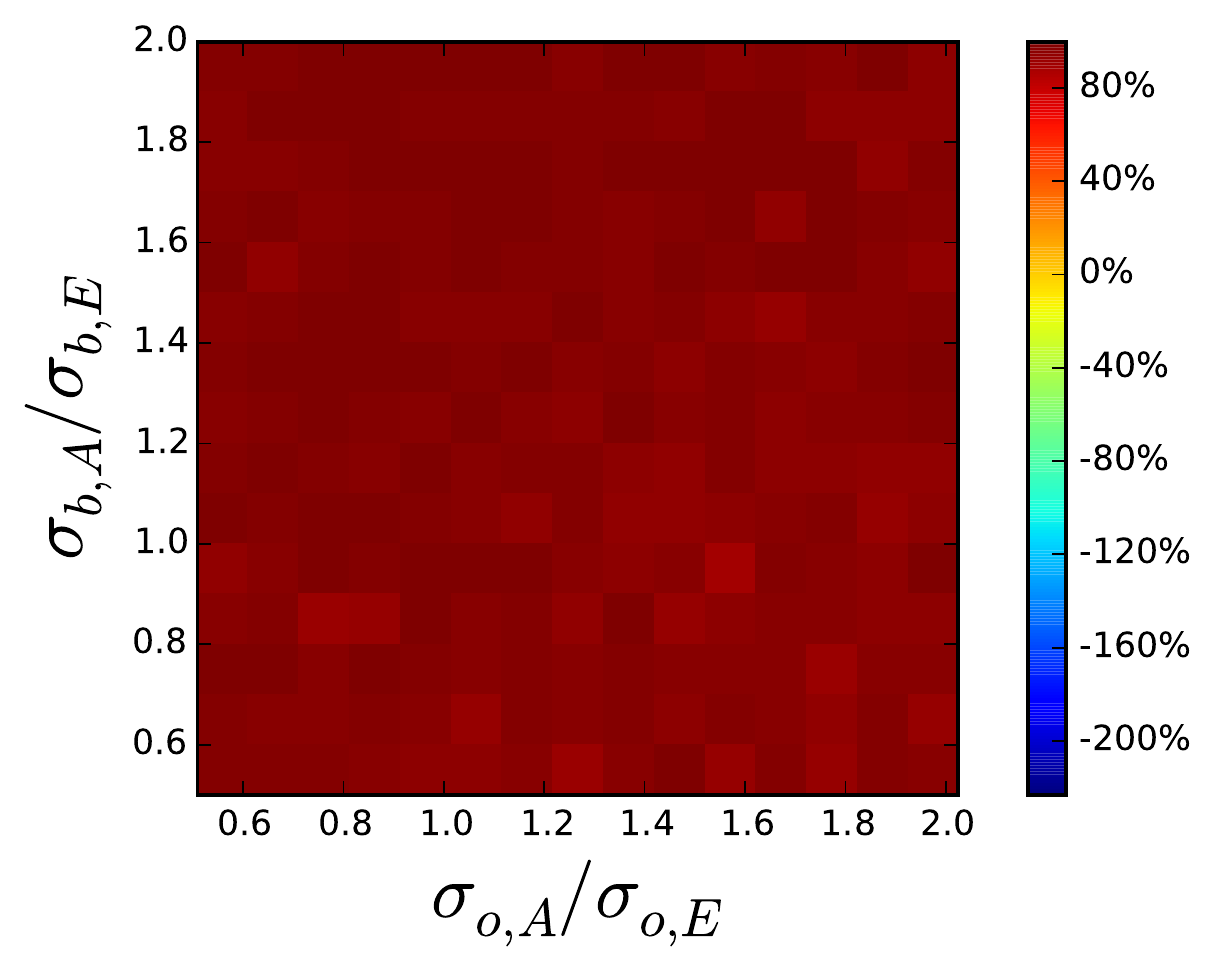}
     \includegraphics[width=2.65in]{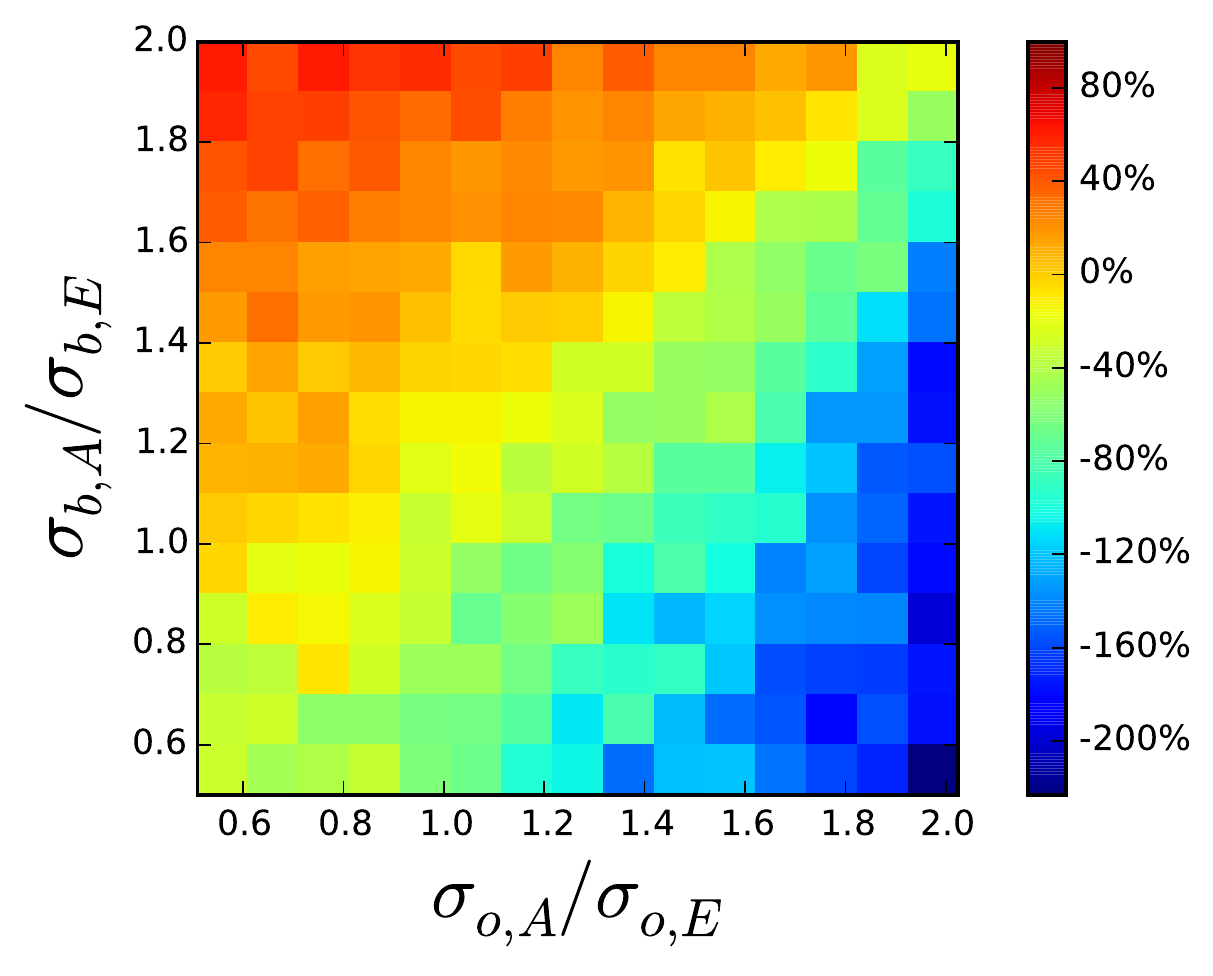}\\
     (a) \hspace{6cm} (b)
     \caption{Average improvement (in \% according to the measures introduced in \ref{sec:Twin experiments set up}) of the background error covariance  $\textbf{B}$  corrected by the proposed localized approach relative to the standard global tuning, (a): with data reduction ($\tilde{\delta_\textbf{B}}$); (b): with data adjustment ($\hat{\delta_\textbf{B}}$); $A$ stands for assumed and $E$ for exact values, respectively, with ($\sigma_{b,E}, \sigma_{o,E} = \sqrt{\sigma_{o,1} \sigma_{o,2}}$) both varying in $[0.025,0.1]$.
     }
  \label{fig:compare_B_under}
\end{figure*}

Fig. \ref{fig:compare_B_under} and Fig. \ref{fig:compare_R_under} collect the results of the Monte-Carlo tests described in \ref{sec:Twin experiments set up} where the ratio of exact error deviation over the assumed one is chosen to vary from $0.5$ to $2$ for both background and observation errors. The improvement in terms of covariance matrices specification is estimated for the standard DI01 algorithm as well as its localized version with two strategies for fitting the observation data. We are interested in the potential advantage of the new methods compared to the standard algorithm. The normalized difference of covariance specification error is drawn as mentioned in \ref{sec:Twin experiments set up} where positive values represent an advantage of localized methods. All tuning methods are applied for $q_\textrm{max}=10$ iterations and we have checked that the sequences $s_{b,q},s_{o,q},s_{b,q}^i, s_{b,q}^i$ have been well converged to 1.

\subsubsection{Measure of improvement of the localized approaches for the estimation  of the background $\textbf{B}$ matrix }

From Fig. \ref{fig:compare_B_under}, one observes that in this test, the localized \di~ with data reduction  always holds a strong advantage (positive value) in terms of matrix $\textbf{B}$ estimation, no matter the exact error deviation, compared to the standard approach. The strategy of data adjustment works well for some parameters combinations, but it becomes less optimal when $\sigma_b$ increases and $\sigma_o$ decreases. Thus careful attention should be brought on the error level of the background state while applying data adjustment strategy. In fact, when the background error level is high, adjustment of the observation data with background state of large variance will take a considerable risk of polluting the observations both in terms of observation accuracy and the knowledge of error covariances.

\subsubsection{Measure of improvement of the localized approaches for the estimation  of the observation $\textbf{R}$ matrix }

\begin{figure*}
  \centering
     \includegraphics[width=2.65in]{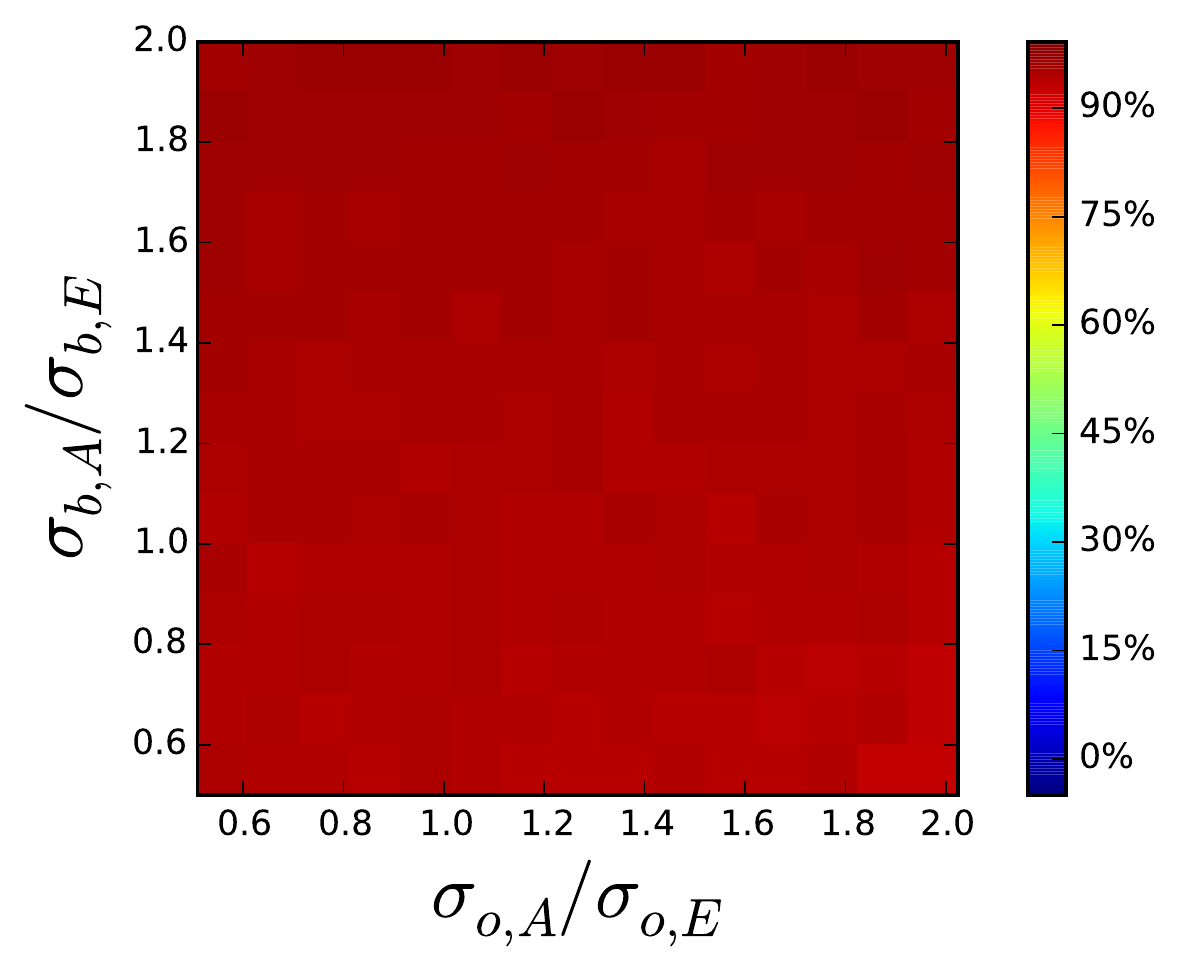}
     \includegraphics[width=2.65in]{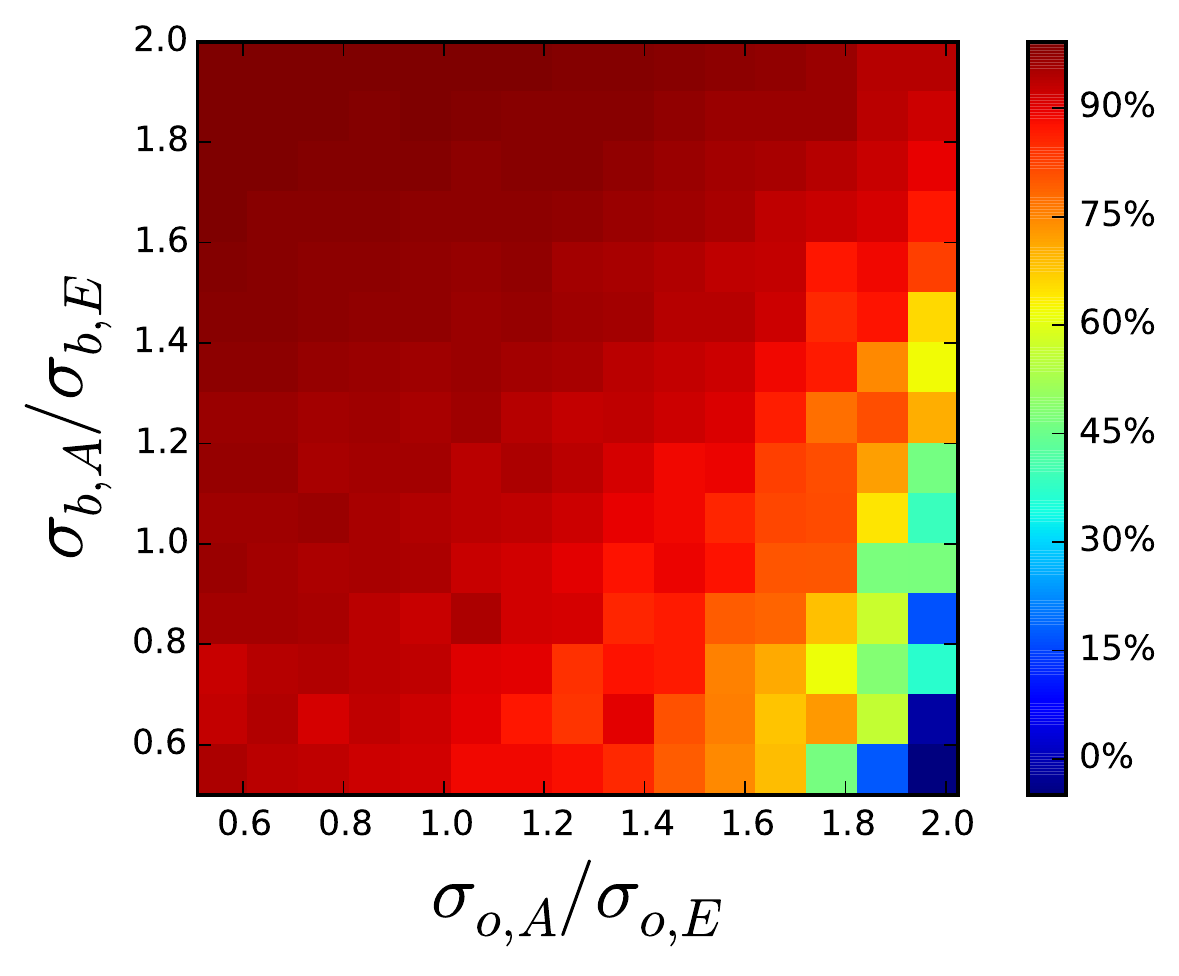}\\
     (a) \hspace{6cm} (b)
          \caption{Same figure as in Fig. \ref{fig:compare_B_under} for observation error covariance improvement $\tilde{\delta_\textbf{R}}$ (a) and $\hat{\delta_\textbf{R}}$ (b).}
  \label{fig:compare_R_under}
\end{figure*}

From Fig. \ref{fig:compare_R_under}, one observes significant advantages in most cases for both new adaptive approaches. In fact, according to the hypothesis of our experiments, the non-homogeneity of observation errors is completely neglected by a standard \di. This non-homogeneity could be covered 
using the graph-based new approach. Similar to the matrix $\textbf{B}$, less optimal results are found when the background error is considerably higher than the observation one.

In these twin experiments, we may conclude that despite the fact that half of the observations are ignored for the covariance tuning, the strategy of data reduction owns in general an  advantage over the one of data adjustment. However, for problems of large dimension, it is possible that most observations are imperfect concerning the correspondence to state communities. Therefore, how to wisely combine these two strategies in real applications for improving the covariance tuning could be a promising topic.

\subsubsection{Test case with a larger difference of error deviation across the two observation clusters}
Similar experiments are also performed with a more significant difference between the two observation groups in terms of their prior error deviation. The setup of experiments is the same as in section \ref{sec:Test case description}, except that the ratio of ${\sigma_{o,E}^{i=1}(\textbf{y}_u)} \Big/{\sigma_{o,E}^{i=2}(\textbf{y}_v)}$ is now set to be 100 instead of 10. 
The same number of experiments as in the previous case are carried out. The test results are summarized in Table \ref{table:2}, according to the cases of under- or over-estimation of prior error amplitude.  As expected, due to the larger difference between the two observation groups and thanks to the assumed  homogeneous  observation matrix $\textbf{R}_\textrm{A}$, the results of the new approaches are even more impressive over a standard DI01 while keeping the same trends against the variations of $\sigma_{b,\textrm{A}} \big/ \sigma_{b,\textrm{E}},  \sigma_{o,\textrm{A}} \big/ \sigma_{o,\textrm{E}}$ similarly to  Fig. \ref{fig:compare_B_under} and Fig. \ref{fig:compare_R_under}.  On the other hand, while the prior estimation of $\sigma_{b,A}, \sigma_{o,A}$ is of extremely poor quality, for example,$\sigma_{o,\textrm{A}} \big/ \sigma_{o,\textrm{E}} > 100 \hspace{3mm} (\textrm{or} < 1/100)$,  we recommend to consider the standard DI01 in the first place.  

\begin{table*}
\centering
\begin{tabular}{ |l||l|l||l|l| }
 \hline
 \bfseries\makecell[c]{Improvement of\\  $\textbf{B}$ and $\textbf{R}$ (in $\%$) } &   \multicolumn{2}{c}{$\sigma_{o,A} < \sigma_{o,E}$}   & \multicolumn{2}{c}{$\sigma_{o,A} > \sigma_{o,E}$} \vline \\
  \hline
 \hline
 observation reduction & \centering $\overline{\gamma_{\tilde{\textbf{B}}}}$& \centering $\overline{\gamma_{\tilde{\textbf{R}}}}$ & \centering $\overline{\gamma_{\tilde{\textbf{B}}}}$ &  $\overline{\gamma_{\tilde{\textbf{R}}}}$ \\
 \hline
 \hline
 $\sigma_{b,A} < \sigma_{b,E}$ & 98.5\% & 96.65\%  & 96.41\% & 96.08\%  \\
 \hline
 $\sigma_{b,A} > \sigma_{b,E}$ & 99.24\% & 96.08\%  & 97.81\% & 96.69\%  \\
 \hline
 \hline
 observation adjustment & $\overline{\gamma_{\hat{\textbf{B}}}}$& $\overline{\gamma_{\hat{\textbf{R}}}}$ & $\overline{\gamma_{\hat{\textbf{B}}}}$& $\overline{\gamma_{\hat{\textbf{R}}}}$ \\
  \hline
 \hline
 $\sigma_{b,A} < \sigma_{b,E}$ & 34.76\% & 99.27\% & -4.03\% & 96.94\%  \\
 \hline
 $\sigma_{b,A} > \sigma_{b,E}$ & 68.22\% & 99.64\% & 30.83\% & 98.81\% \\
 \hline
\end{tabular}
\caption{Averaged gain improvement of error covariances ($(\textbf{B}, \textbf{R})$in $\%$) with ${\sigma_{o,E}^{i=1}(\textbf{y}_u)} = 100 {\sigma_{o,E}^{i=2}(\textbf{y}_v)}$ via two graph clustering localization strategies (observation reduction and observation adjustment).  Both  $\sigma_{b,E} $ and $\sigma_{o,E}$ vary in $[0.025,0.1]$.}
\label{table:2}
\end{table*}

\section{Discussion}
\label{sec:discussion}

Localisation technique is an important numerical tool which contributes to the success of solving high-dimensional data assimilation problems for which ensemble estimates are unreliable. It is based on the assumption that correlations between dynamical system variables eventually decay with the physical distance. This simple rationale is put to use either to make the assimilation of observations more local  (domain localization) or to numerically impose a tapering of distant spurious correlations (covariance localization) and leads to very different implementations and numerical difficulties. Domain localization is interesting because it makes the problem more scalable and the implementation more flexible in the sense that the original global formulation can be broken-up into several smaller subproblems. Nevertheless, the assimilation of {\em non-local} observations and/or observations from different sources and at different scales becomes increasingly important for instance due to the use of satellite observations.\\ In this work, we propose to generalize the concept of domain localization relying on graph clustering state decomposition techniques. The idea is to automatically detect and segregate the state and observation
variables in an optimal number of clusters, more amenable to scalable data assimilation, and use this decomposition to perform efficient adaptive error covariances tuning cluster-wise.
This unsupervised localization technique based on a {\em linearized} state-observation measure is general and does not rely on
any prior information such as relevant spatial scales, empirical cut-off radius or
homogeneity assumptions. In this paper, the Fluid method is chosen for
applications because of its computational simplicity, especially for sparse graphs.
In terms of covariance diagnosis, the \di~ is chosen because available data is often limited for our industrial applications. \\
The methodology is applied to a simple twin experiments data assimilation problem for which the Jacobian matrix of the observation operator is chosen to reflect a dual clustering of the state-observation mapping; the components of which are then randomly mixed in order to hide any particular structural pattern. Simply speaking, there exist two {\em hidden} communities of state variables, each of them preferably connected to their own observations community. The problem is far from trivial as — the segregation resulting in clustering is not related to let us say spatial separations, — exact background error magnitude is supposed to be homogeneous in our tests but the clusters have different exact observation errors and also because — there exists some inter-connectivity between the clusters. Considering the latter, two simple numerical approaches are proposed in order to handle a data {\em reduction} or a data {\em adjustment} strategy. The problem is investigated for a wide range of assumed prior covariances and the graph clustering approach with adaptive covariance tuning is much more efficient than a global adaptive covariance tuning approach, especially in the case of \di~.\\
Here follows some discussion about the method. The graph clustering algorithm uses an adjacency matrix derived from a linearization of the observation operator. Therefore, it seems reasonable to anticipate that the approach will be more appropriate for linear or weakly nonlinear problems. For time-dependent strongly nonlinear problems, one may need to rely on the community detection algorithm multiple times, which could be computationally expensive.  \\
Another critical point relates to the inter-cluster connectivity which materialize the fact that real applications problems will never be fully separable. Here, we have made the choice to circumvent the difficulty by disposing of the troublesome shared observations. Nevertheless, this approach might be impractical for applications with a large number of clusters and overlaps. In this case, alternate strategies will have to be considered, much likely involving a search for optimal ordering of the subspaces covariance tuning.\\   
Finally, our localization approach will perform better if the assimilation problem, represented by a graph, is well separable under our cluster analysis; i.e. in the sense that the data assimilation problem is decomposed into a certain number of subsets problems minimizing the overlap between subsets.  This will somewhat depend on the graph cluster analysis algorithm retained but more predominantly on the chosen measure of similarity. For the former, it will be useful to monitor some performance metrics as a function of the number of clusters for a given graph clustering algorithm. In this work, we base the measure of similarity on the linearized observation operator. Complementary to this approach, it may be interesting to combine a measure involving prior knowledge of error covariances with the state-observation mapping, i.e. $\left | \textbf{H}\right | \textbf{B} \left |\textbf{H} \right |^T$ instead of $\left | \textbf{H}\right | \left |\textbf{H} \right |^T$. 
This might provide a way to scalable optimization of covariance structures between observations and model variables instead of covariance structures in the prior alone. After this methodological contribution, future work will consider applying these methods to more challenging real industrial applications.
\bibliographystyle{abbrv}
\bibliography{these_siboCHENG}

\end{document}